\documentclass[nospthms]{svjour3}
\usepackage{graphicx}
\usepackage{amssymb}
\usepackage{amsmath,amscd}

\usepackage{url,verbatim}
\usepackage{pgfplots}
\usepackage{caption,subcaption}
\usepackage[nottoc,notlot,notlof]{tocbibind}

\RequirePackage[colorlinks,citecolor=blue,urlcolor=blue]{hyperref}
\RequirePackage[colorlinks,citecolor=blue,urlcolor=blue]{hyperref}
\usepackage{breakurl}

\newtheorem{theorem}{Theorem}

\newcommand\EE{{\mathbb E}}

\newcommand\HH{{\mathbb H}}
\newcommand\RR{{\mathbb R}}
\newcommand\ZZ{{\mathbb Z}}

\newcommand\sD{{\mathcal D}}

\newcommand\var{\text{\rm var}}

\newcommand\NN{{\mathbb N}}

\newcommand\PP{{\mathbb P}}

\newcommand\TT{{\mathbb T}}

\newcommand\pc{p_{\text{\rm c}}}
\newcommand\pslab{p_{\text{\rm slab}}}
\renewcommand\o{\text{\rm o}}
\newcommand\q{\quad}
\newcommand\qq{\qquad}
\newcommand\rad{\text{\rm rad}}
\newcommand\lra{\leftrightarrow}
\newcommand\wt{\widetilde}
\newcommand\om{\omega}
\newcommand\eqd{\mathrel{\stackrel{\text{\rm d}}{=}}}
\def\ol{\overline}
\def\what{\widehat}
\newcommand\Bs{B_{\text{\rm s}}}
\newcommand\sP{{\mathcal P}}
\newcommand\sL{{\mathcal L}}
\newcommand\es{\varnothing}

\newcommand\oo{\infty}
\newcommand\resp{respectively}
\newcommand\lest{\le_{\text{\rm st}}}

\newcounter{mycount1}\newcounter{mycount2}\newcounter{mycount3}\newcounter{mycount}

\newenvironment{numlist}{\begin{list}{\rm\arabic{mycount2}.}%
   {\usecounter{mycount2}\labelwidth=1cm\itemsep 0pt}}{\end{list}}
\newenvironment{letlist}{\begin{list}{\rm(\alph{mycount3})}%
   {\usecounter{mycount3}
   \labelwidth=20pt
   \itemindent=-2pt
   \labelsep=5pt
   \leftmargin=1cm
   \itemsep 0pt}}{\end{list}}

\numberwithin{equation}{section}
\numberwithin{theorem}{section}
\numberwithin{figure}{section}

\begin{document}
\title{Harry Kesten's work in probability theory}

\author{Geoffrey R.\ Grimmett}
\institute{Statistical Laboratory, Centre for
Mathematical Sciences, Cambridge University, Wilberforce Road,
Cambridge CB3 0WB, UK; School of Mathematics \&\ Statistics,
The University of Melbourne, Australia;
Heilbronn Institute for Mathematical Research, Bristol University, UK.
\email{grg@statslab.cam.ac.uk} }

\dedication{In memory of Harry Kesten, inspiring colleague, valued friend}


\begin{abstract}
We survey the published work of Harry Kesten in probability theory,
with emphasis on his contributions to
random walks, branching processes, percolation, and related topics.  
A complete bibliography is included of his publications. 
\keywords {Probability, random walk, branching process, random matrix, 
diffusion limited aggregation, percolation.} 
\subclass {60-03, 60G50, 60J80, 60B20, 60K35, 82B20.}
\end{abstract}

\date{Submitted 20 April 2020, revised 26 October 2020}

\maketitle

\tableofcontents

\section{Overview}\label{sec:1}

Harry Kesten was a prominent mathematician and personality 
in a golden period of probability theory from 1956 to 2018. 
At the time of Harry's move from the Netherlands to the USA in 1956, 
as a graduate student aged 24, much of the foundational
infrastructure of probability was in place. 
The central characters of probability had long been identified
(including random walk, Brownian motion, the branching process,
and the Poisson process), 
and connections had been made and developed between
\lq pure theory' and cognate areas ranging from physics to finance.
In the half-century or so since 1956,
a coordinated and refined theory has been developed, and probability has been recognised
as a crossroads discipline in mathematical science. Few mathematicians have contributed as much during this 
period as Harry Kesten.

Following a turbulent childhood (see \cite{HK-AMS}), Harry studied mathematics with 
David van Dantzig and Jan Hemelrijk in Amsterdam, where
in 1955 he attended a lecture by Mark Kac entitled \lq\lq Some probabilistic aspects of potential theory".
This encounter appears to have had a decisive effect, in that Harry moved in 1956 to Cornell University
to work with Kac. In due course, and together with his colleagues including 
his fellow emigr\'es Frank Spitzer
and later Eugene Dynkin, Harry's work and influence supported the Cornell Mathematics Department as a 
leading institution worldwide in probability and beyond. 
Many visitors  were attracted
to this extraordinary academic niche in upstate New York, where they were received with warmth,
and invited to participate in a variety of mathematical and physical activities.
   
Mark Kac had a lasting influence on Harry, and it was apposite
that Harry should write an appreciation of the former's work on his death in 1984 (as he did for Spitzer
in 1993 (see \cite{MR866337,MR1217555}\footnote{Citations beginning with the 
letter K refer to the list of Harry's publications at the end of this article.}).
Coincidentally, Mark Kac's first publication was dated the year of Harry's birth, 1931, and
his last appeared shortly after one of Harry's best known results,
namely his proof in 1980 that the critical probability of bond percolation on $\ZZ^2$
equals $\frac12$.

A number of further individuals influenced Harry's work, including 
particularly his friend and colleague Frank Spitzer,
and more distantly John Hammersley, whose earlier papers on self-avoiding walks,
percolation, first-passage percolation, and subadditivity are frequently reflected in Harry's own papers.
Harry had in common with Hammersley  a preference for what the latter termed 
\lq\lq implicated'' mathematics over the \lq\lq contemplative'' sort, and the
solving of problems featured strongly in both their scientific lives (see \cite[p.\ 1132]{gw}).

While Harry's publications lie largely within probability theory, his broader expertise 
extended into neighbouring areas of mathematics, including topics in 
algebra, analysis, geometry, and statistical physics. 
Already by 1960 he had substantial results for random walks on groups,
products of random matrices, and aspects of what we call here
\lq probabilistic Diophantine approximation'. 
These interests expanded over the next decade
to include random walk and potential theory, stable and L\'evy processes, self-avoiding walks,
and branching processes, and in the 1980s  to random walk in random environments,
and random recurrence relations. 

Harry's career was punctuated in 1980 by his proof that $\pc=\frac12$ for 
bond percolation on the square lattice.
Before then, he worked largely on random walks and branching processes, and afterwards 
his interests shifted towards mathematical processes inspired by physics. 
Common themes and strands are discernible. For example, 
subadditivity crops up first in his 1960 paper on random matrices, and much later in first-passage 
percolation and random flows, with self-avoiding walks in between.   
Local surgery and exponential estimates are techniques which he used in a number of contexts.

He had a taste for hard problems, 
and he developed a fearsome reputation as a problem-solver.
He would plough on beyond the capacities of others when confronted with technical difficulties or complications,
and it could be years before the community caught up with him. 
A number of his major works have later been 
overtaken by \lq neater' proofs of others, though his original methods, being more \lq hands on',
can prove more robust as the assumptions are perturbed. 

Certain decisions had to be confronted in writing this summary of Harry's work.
On the one hand, a chronological account allows the reader to witness the development of Harry's interests, and 
more generally the evolution of probability theory; on the other hand, many readers may find his work on physical systems 
more approachable than his somewhat more technical
work on the mathematics of random walks and the like. The strategy adopted in
the current article is to present the author's choice
of highlights of
Harry's research (in Section \ref{sec:top10}) 
followed by the fuller, more chronological survey
of Sections \ref{sec:rw}--\ref{sec:perc}.
  
The time-ordered framework of this article is subject to frequent local deviations. 
Thus, we begin in Section \ref{sec:rw} with random walks in their various forms,
followed by products of random matrices, self-avoiding walks, and branching processes;
thence to classical and first-passage percolation, and the many related problems
that attracted Harry. Certain later results are dovetailed into earlier sections,
where they risk disturbing the chronological flow with brief reviews.
Some important papers that do not fit easily into this programme
are deferred to the final Section \ref{sec:last}. 

The author acknowledges the incompleteness of this personal perspective, and he pleads
limited knowledge and lack of space. He apologises to those whose
relevant work is not listed. Special mention is made of Rick Durrett's excellent article \cite{MR1703122}
written for the 1999 Festschrift for Harry (see \cite{hk-fest}),
which has been a support in the writing of the current work.

Not quite all of Harry's work is publically available. In addition to his 
two as yet unpublished articles on \url{arxiv.org},
which are  included in the appended bibliography, there is a handful of papers whose lengths he 
considered to be disproportionate to their novelty. He kept these in his bottom drawer,
and we have abided by his judgement. 

Mention is made of some of Harry's professional activities beyond his own research. 
In addition to serving on the Editorial Boards of a number of journals
(including the current journal, then known as ZfW), he was the Managing Editor of The
Annals of Probability from 1982--1984, in fact the fourth such Editor since the journal's inception in 1972.
He served as an Elected Member of the Council of the Institute for Mathematical Statistics from 1988--1991.

Harry Kesten's impact on mathematics was broader than 
his decidedly outstanding written output. He was in addition a very popular
and respected figure, and an excellent collaborator and correspondent.
As a role model and friend for younger colleagues, his positive influence on the intellectual
and personal values of the probability community is amongst his greatest legacies.  

\section{Highlights of Harry Kesten's research}\label{sec:top10}

This section contains very brief accounts of some of the highlights of Harry's work,
taken out of chronological order and context. Readers may follow the links to
the fuller accounts presented in Sections \ref{sec:rw}--\ref{sec:last}.

\subsection{Percolation}

\subsubsection*{Critical probability of bond percolation on $\ZZ^2$ 
{\rm(\S \ref{sec:perc})}}
Harry's proof that $\pc=\frac12$ for bond percolation on $\ZZ^2$ was a watershed for the subject.
It solved a notorious old problem and introduced a number of new ideas.
In so doing, it illuminated percolation theory as a key process of probability,
and it contributed to a development of disordered systems that is very active at the time of writing.

\subsubsection*{The van den Berg--Kesten inequality {\rm(\S \ref{sec:bk})}}
Correlation inequalities are key to the study of interacting systems. The so-called
FKG inequality has been especially important, 
and it asserts that increasing events are positively correlated
(for suitable probability measures). \emph{Negative correlation} is a more elusive concept
than its positive sibling. Probably the most important negative correlation inequality for product measures is
that of van den Berg and Kesten from 1985. They introduced the notable concept of 
the \lq disjoint occurrence' of events, 
and they proved the associated inequality.  It has been very useful since.

\subsubsection*{Uniqueness of the infinite cluster {\rm(\S \ref{sec:uic})}}

How many infinite open clusters does a supercritical percolation process 
on $\ZZ^d$ possess? Aizenman, Kesten, and Newman showed the 
answer to be a.s.\ one. They introduced methods that have since been useful
elsewhere, and they opened the door to a series of important papers 
by others on disordered systems on graphs of a variety of types.
It turns out that the answer is always one for amenable graphs,  whereas
the picture is more diverse for non-amenable graphs.

\subsubsection*{Incipience and subdiffusivity {\rm(\S \ref{sec:incic})}}

Critical percolation is believed (and in most dimensions known) to have only finite open clusters.
From physics came the idea of studying the infinite percolation cluster that critical percolation is \lq trying to form'.
Harry made sense of this notion in two dimensions with 
his proof of the existence of a probability measure defined as the limit
of the (critical) product measure conditioned on the origin lying in an increasingly large cluster. He called
this limit the \emph{incipient percolation cluster} $\wt C$.

He was able to show a fractal-like structure of $\wt C$, in part
by utilising de Gennes' proposal of setting in motion a random walker. 
Whereas the natural normalization of a random walk $(X_n)$ 
on $\ZZ^2$ is $\sqrt n$, it turns out that, for a walk $(X_n)$ on $\wt C$, there exists
$\epsilon>0$ such that $n^{-\frac12+\epsilon}X_n\to 0$. This \lq subdiffusivity' occurs
because the walk spends a lot of its time in blank alleys of $\wt C$.

\subsubsection*{Scaling relations in two dimensions {\rm(\S \ref{sec:scal})}}

There is a full physical picture of the percolation phase transition in two dimensions. A great deal
is \lq known' but, even today, relatively little has been proved (with the 
exception of site percolation on the triangular lattice). It is expected that functions have 
power-law singularities at and near the critical point $\pc$, and that their \lq critical exponents' satisfy the
so-called scaling relations.

A major step forward was taken by Harry in 1987 towards an understanding of such exponents and relations.
Under the assumption of the existence of two exponents, he proved the existence of  five more, and the validity
of the associated scaling relations. On the way to this, he studied the so-called \lq arm events' that
have been so important in more recent work.
This article has had a major impact on the rigorous theory
of the percolation phase transition in two dimensions.

\subsubsection*{Random flows {\rm(\S \ref{sec:fpp})}}

\lq First-passage percolation' is the study of the rate of spread of material through a random medium. 
Harry made a systematic study of the related problem of the maximal flow through a medium subject to
random edge-capacities.  Careful control of
dual surfaces is required in three or more dimensions.

\subsection{Random walks}

\subsubsection*{Random walks on groups {\rm (\S \ref{sec:rwgp})}}

In his distinctively innovative PhD thesis from 1958, Harry introduced the theory of random walks on countable discrete groups.
He proved that the probability of return to the starting point after $2n$ steps behaves in the manner of
$\lambda^{2n}$ where $\lambda$ is the spectral radius. Moreover, $\lambda=1$ if and only if
the group is amenable. This condition for amenability has become one of the tools of
geometric group theory, and the field that he initiated remains active today.

\subsubsection*{Random walk in random environment {\rm(\S \ref{sec:rwre})}}

Undergraduates learn about random walks on $\ZZ$ that, when at $x\in\ZZ$,  jump one step to the right (\resp, left) 
with probability $\alpha$ (\resp, $\beta=1-\alpha$). Life becomes much more complicated if 
the value of $\alpha$ at the point $x$ 
is a random variable $\alpha_x$. It turns out that the long-term
behaviour of the walker is heavily influenced by domains $D$ of $\ZZ$ where the $(\alpha_x: x \in D)$ conspire
to slow it down.   With Kozlov and Spitzer, Harry proved a very precise result about
the walker's position, in terms of a stable law of index $\kappa$,
where $\kappa$ is given by $\EE((\beta_0/\alpha_0)^\kappa)=1$.

\subsubsection*{Diffusion limited aggregation {\rm(\S \ref{sec:dla})}}

This famous model for aggregation has resisted many attempts to build rigorous theory. 
Drifting particles aggregate at the first point of a composite that they hit. 
Simulations support the belief that growth is dendritic and fractal. In 1987 Harry proved 
one of the few rigorous results, namely that (in two dimensions)
the radius of the composite grows no faster than $n^{\frac23}$.

\subsection{Branching processes and $L\log L$ (\S \ref{sec:bp})}
Harry liked necessary and sufficient conditions, as exemplified in his most prominent works on branching processes.
Consider a supercritical branching process $(Z_n)$ with one progenitor and mean family-size
$\mu\ge 1$. The limit $W=\lim_{n\to\oo} Z_n/\mu^n$ exists, by martingale convergence. With Stigum, Harry found the
necessary and sufficient condition for $\EE(W)=1$, namely $\EE(Z_1\log^+ Z_1)<\oo$.
Around the same time, with Ney and Spitzer, he established necessary and sufficient 
conditions for the Yaglom and Kolmogorov
laws for critical branching processes.

\subsection{Products of random matrices (and scalars) (\S \ref{sec:prm})}

In early work with Furstenberg from 1960, Harry investigated the limiting behaviour of
the partial products $Y_n$ of a stationary, ergodic sequence $(M_n)$ of random matrices.
They proved in particular that the limit of $n^{-1}\log \|Y_n\|$ exists, 
and furthermore that a central limit theorem holds under suitable assumptions.

Harry studied later the stochastic  recurrence equation $Y_n=M_nY_{n-1}+Q_n$, 
as well as the apparently simpler case of a \emph{scalar} random variable $Y$ such that $Y$ and $MY+Q$ have the
same distributions. His proof that $Y$ is generically heavy-tailed has had enormous influence
in numerous areas of applied probability, statistics, and finance.

\subsection{Self-avoiding walks (\S \ref{sec:saw})}

What can be said about the number $\chi_n$ of $n$-step self-avoiding walks (SAWs) on a lattice
$\sL$, with a given starting point?
The combinatorics of SAWs are especially complicated, and asymptotics for $\chi_n$ are challenging to derive.
It is classical that $\chi_n^{1/n}\to \kappa$ for some \emph{connective constant} $\kappa$ associated with the
lattice $\sL$. When $\sL=\ZZ^d$,
Harry proved the more refined  ratio limit theorem $\chi_{n+2}/\chi_n \to \kappa^2$, and in so doing 
introduced a \lq pattern theorem' which has had lasting impact in this challenging and visible field.

\subsection{Polar points of L\'evy processes (\S \ref{sec:gtrw})}

Let $(X_t)$ be a $d$-dimensional L\'evy process, and define the hitting probability
$h(r)=\PP(X_t=r\text{ for some }t>0)$ for $r \in \RR^d$. The point $r$ is called \emph{polar}
if $h(r)=0$. The identification of the polar points became a prominent problem in the 1960s.
Harry identified the set of polar points in one dimension, and proved 
(subject to reasonable conditions) that polarity is the norm in
two or more dimensions.

\subsection{Probabilistic Diophantine approximation (\S \ref{sec:eucl})}

Perhaps inspired by an early collaboration with Mark Kac, Harry wrote several papers around 1960
on aspects of Diophantine approximation, using frequently the language and methods of probability theory. 
His best known work is a  paper which set the scene for the more recent theory of   
bounded remainder sets. Let $\{x\}=x-\lfloor x\rfloor$ denote the fractional part of the real number $x$.
With $0\le a< b\le 1$ and $b-a<1$, let $S_n(x)$ be the number of integers $k\in\{1,2,\dots,n\}$ such that $\{kx\} \in[a,b)$.
In proving that $S_n(x)-n(b-a)$ is bounded in $n$ if and only if
$b-a=\{jx\}$ for some integer $j$, Harry answered a question of Erd\H os and Sz\"usz.

\section{Random walk}\label{sec:rw}

\subsection{Random walk}

The purist might view probability theory as finite measure theory plus conditional
expectation. The concept of \lq independence' is fundamental to the theory,
as is its negation \lq dependence'. As Kac wrote in \cite{kac-carus}:
\lq\lq This notion [statistical independence] originated in probability theory and for a
long time was handled with vagueness which bred suspicion
as to its being a bona fide mathematical notion.'' Independence was considered historically
as a basis for the multiplication of probabilities, and it seems to have been Steinhaus (student of Hilbert, 
advisor to Banach and Kac) who introduced
the now familiar definition of independence for a countable family of random variables (as reported by Kac \cite{kacI};
see also  \cite[p.\ 1104]{MR866337}).

Random walk is the theory of the partial sums of a 
sequence of independent, identically distributed  (iid) random variables.
Such partial sums are viewed as a process, suitably indexed, and the emphasis is generally upon the geometrical
aspects of this process. Random walk is arguably the most fundamental process of probability.
Recent books on the topic include \cite{lawler-int,lawlim}.

\subsection{Random walk on groups}\label{sec:rwgp}

It is not surprising, given the circumstances of his arrival at Cornell, that Harry's PhD thesis
was concerned with random walk, albeit in the novel context of groups. In his 35 page thesis of 1958
(\cite{MR2612793}, published as \cite{MR109367}),
he initiated the theory of random walks on groups, a healthy topic even 60 years later. 

The usual random walk takes jumps around the abelian groups $\ZZ^d$ or $\RR^d$ 
under addition. More generally, let $G$ be a countable
(additive) group, and let $A=\{a_1,a_2,\dots\}$ be a generating set for $G$ and $P=(p_1,p_2,\dots)$
a strictly positive vector with sum $\frac12$. The triple $(G,A,P)$
gives rise to a random walk $X=(X_n: n \ge 0)$ 
on $G$ that, when at $g \in G$, moves at the next step to
$g \pm a_i$ with probability $p_i$ (for each $i$). Harry investigated the relationship 
between the spectrum of the transition matrix $M$ (considered as an operator on $l^2(G)$)
and the properties of the group $G$. The \emph{spectral radius} $\lambda(G,A,P)$ of $M$ is
the supremum of $|\lambda|$ taken over all eigenvalues $\lambda$ of $M$. 
Write $0$ for the zero element of $G$, and $\PP_x$ for the law of $X$ given $X_0=x$. 

\begin{theorem}  {\bf\cite{MR109367}}
\begin{letlist}
\item The spectral radius $\lambda(G,A,P)$ equals the maximal eigenvalue of $M$.
\item We have that
$\lambda(G,A,P) = \lim_{n\to\oo} \PP_0(X_{2n}=0)^{1/(2n)}$.
\item
Whether or not $\lambda(G,A,P)=1$ depends only on the 
structure of the group $G$, in the following sense.
We have that $\lambda(G,A,P)=1$ if and only if
$\lambda(G,B,Q)=1$ for some generating set $B$ and strictly positive $Q$.
\end{letlist}
\end{theorem}

Write $\lambda(G)=1$ when the common value of part (c) satisfies $\lambda(G,A,P)=1$.
The following criterion for amenability is a fundamental result in geometric group theory.

\begin{theorem} {\bf\cite{MR112053}}
The group $G$ is amenable if and only if $\lambda(G)=1$.
\end{theorem}

In the collaboration that followed, Kesten and Spitzer \cite{MR195163} 
asked about potential kernels for random walks
on countably infinite abelian groups. They presented a criterion for recurrence, and they established a potential theory.
In his Berkeley symposium paper \cite{MR0214137}, Harry collected together a number of related open
problems, including what came to be known as \lq Kesten's conjecture': the only recurrent groups are the finite
extensions of $\{0\}$, $\ZZ$, and $\ZZ^2$. This prominent conjecture was finally verified in its generality in 1986 by
Varopoulos \cite{varop} (see \cite{guiv}).

Random walks on groups have attracted much attention in recent years; see \cite{woess}. 
In the case of finite groups,
it is appropriate to mention Saloff-Coste's chapter \cite{saloff}, included
in the collection \cite{MR2023649} edited by Harry.

\subsection{General theory of random walk}\label{sec:gtrw}

The decade of the 1960s was a rich period for the theory of random walk, with the first edition of Spitzer's important book  \cite{Spitzer-RW}
published in 1964. Harry wrote a number of significant articles around this time, 
of which four were concerned with ratio limit theorems.
 
\begin{theorem} 
Let $X$ be an irreducible random walk on $\ZZ^d$ where $d \ge 1$. Assume $X_0=0$,
and let $T$ be the time of the first return by $X$ to $0$. 
\begin{letlist} 
\item {\bf\cite{MR162279}} The probabilities $r_n=\PP_0(T>n)$
satisfy $r_{n+1}/r_n \to 1$ as $n\to\oo$.

\item  {\bf\cite{MR142160}}
Furthermore, the limit
$$
\lim_{n\to\oo} \frac{\PP_x(T > n)}{\PP_0(T > n)}=\delta_{0,x} + a(x)
$$
exists for $x \in \ZZ^d$, where $a(x)$ is  the potential kernel
$$
a(x) = \sum_{n=0}^\oo \bigl[\PP_0(X_n=0)-\PP_x(X_n=0)\bigr] < \oo.
$$
\end{letlist}
\end{theorem}

Part (a) is a key part of the proof of part (b), due to Kesten, Ornstein, and Spitzer.
The main purpose of \cite{MR162279} was to extend part (b) to the number of visits
of $X$ to an arbitrary finite subset $\Omega$ of $\ZZ^d$, that is, to exhibit an expression
for the limit of the ratio $\PP_x( n,r,\Omega)/\PP_0(T>n)$ as $n\to\oo$, where the numerator is the probability
that $X$ makes exactly $r$ visits to $\Omega$ up to time $n$. This is a ratio of sums, and 
one may in principle obtain more refined results by examining the ratios of certain individual terms in these sums.
This last project was pursued by Harry alone in \cite{MR163365}. It was left there as a conjecture
to show that $f_n^{(r)}:= \PP_0(n,r,\{0\})-\PP_0(n-1,r,\{0\})$ satisfies
$$
\lim_{n\to\oo} \frac{f_n^{(r)}}{f_n^{(1)}} = r,\qq r=2,3,\dots
$$
the proof of which (for aperiodic, recurrent, symmetric chains) was delayed until \cite{MR270454}.
 
While at The Hebrew University in 1960, Harry began his work on
$\alpha$-stable processes $X$, in an investigation of the intervals of time during which
$X_t>0$.  He derived asymptotics for the number $N(\epsilon)$ of such intervals with length greater than $\epsilon$, 
and observed that his answers for $\alpha\ne 1$ coincided with results of Chung and Kac 
\cite{chungk1} on 
changes of sign in sums of iid variables with a stable distribution. This led to a correspondence with Chung and Kac about the validity
of their claims when $\alpha=1$ (see \cite[p.\ 1119]{MR866337}). Harry wrote up his results in \cite{MR0148096}.

In the two linked papers \cite{MR0124934,MR0125639}, Harry extended earlier results of Spitzer and Stone \cite{spst}
concerning the mean number of visits by a random walk on $\ZZ$  to a given state before
exiting a given bounded interval. He was principally motivated to consider random walks with
$\alpha$-stable jumps. This work led naturally to two subsequent papers, the first
of which, \cite{MR178494}, answered a question of Erd\H os with an iterated logarithm law
for the number of visits by a random walk to its most visited state.
In the later paper \cite{MR234509}, 
Harry proved that, if the mean number of visits to a given interval $I$ up to time $n$
grows in the manner of $n^{1-1/\alpha}$ with $\alpha\in(1,2]$, then the jumps belong 
to the domain of attraction of a symmetric $\alpha$-stable law.  This answered a question of Spitzer
concerning a converse result of Darling and Kac \cite{dkac}. 

Let $X$ be a $d$-dimensional L\'evy process, that is, a random process taking values 
in $\RR^d$ with stationary independent increments.
Announced in \cite{MR251797} and proved in his AMS Memoir \cite{MR0272059} are Harry's
results on the question: when is it the case that the hitting probability 
$h(r) = \PP(X_t=r \text{ for some } t>0)$
satisfies $h(r)>0$?  He identified seven cases when $d=1$, and included a criterion  
(see Theorem \ref{thm:polar}) for the statement $h(r)=0$ for all $r\in \RR$,
while warning the reader that \lq\lq In most practical situations one will have a hard time applying the criterion $\dots$''. 
The cases $d \ge 2$ are more complicated still. His $d=1$ result 
was motivated by and answered a question of Chung on solutions to a certain convolution equation.

From the spread of results in his paper, we extract one in the following paraphrase. 
The point $r\in\RR^d$ is called \emph{polar} if
$h(r)=0$, and we let $\sP$ denote the set of polar points.

\begin{theorem} {\bf\cite{MR0272059}}\label{thm:polar}
Let $X$ be a one-dimensional L\'evy process which is not a compound 
Poisson process. Then $\sP=\RR$ if and only if the characteristic exponent $\Psi$ of the process $X$ satisfies
$$
\int_{-\oo}^\oo \Re\left(\frac1{1+\Psi(u)}\right)\,du=\oo,
$$
where $\Re(z)$ denotes the real part of $z$.
\end{theorem} 

Shortly after the publication of this paper,
Bretagnolle \lq\lq obtained very elegant and powerful probabilistic arguments 
and thereby obtain[ed] considerably simpler proofs of the 
[Harry's] results'' (quotation from Harry's review of Bretagnolle's paper \cite{Bret}).
Of the numerous modern accounts of L\'evy processes, we mention the books \cite{bertoin,doney,kyp}. 

Harry delivered the 1971 Rietz Lecture of the Institute of Mathematical Statistics.
In the resulting article \cite{MR301786}, he surveyed three aspects of the theory 
of the sum $S_n$ of $n$ independent random variables $X_i$, namely, (i) spread of distributions and concentration functions,
(ii) ratio limit theorems and potential kernels, (iii) the
set of accumulation points of $S_n/\gamma_n$ for suitable $\gamma=(\gamma_n)$, when the $X_i$ are iid.
His inclusion in the title of \lq\lq without moment conditions'' (and omission of \lq\lq identically distributed''), 
is typical of his inveterate search for generality, but he resisted a
temptation to move much beyond one dimension.
Parts (i) and (ii) are important surveys of earlier distributional theory of Harry and others,
including \cite{MR258095} and the material on ratio limit theorems remarked above. 

Part (iii) of \cite{MR301786} is concerned with the \lq strong' theory of $S_n/\gamma_n$.  
Let 
$$
A(S,\gamma) =   \bigcap_{m=1}^\oo \overline{\bigl\{S_n/\gamma_n: n \ge m\bigr\}}.
$$
For any given distribution function $F$ of the $X_i$, there exists a 
deterministic closed subset $B(F,\gamma)$ of
the extended real line $\ol\RR$ such that $A(S,\gamma)=B(F,\gamma)$ a.s.
The strong law of large numbers
and the law of the iterated logarithm may be viewed 
as determinations of $B(F,\gamma)$ for $\gamma_n=n$ and $\gamma_n=\sqrt{2n\log\log n}$
under suitable conditions.  Harry ended his paper with a very general law of the iterated logarithm
for the sums of an iid sequence, subject to the necessary and sufficient condition that 
their common distribution $F$ lies in the domain
of partial attraction of the normal law. 

\begin{theorem}\label{sl1}
{\bf\cite{MR266315,MR301786}}
\begin{letlist}
\item If $B(F,(n))$ contains two or more points, then it contains $+\oo$ and $-\oo$.
\item For any closed subset $B$ of $\ol\RR$ containing $\{-\oo,+\oo\}$, there exists $F$ such that $B(F,(n))=B$.
\end{letlist}
\end{theorem}

Harry asked about the set $B(F,(n^\alpha))$, and he returned
to this question in joint work with Erickson.
Let $\Bs(F,(n^\alpha))$ be the set of $b\in\ol\RR$ such that there exists a non-random sequence $n_k$ 
satisfying $S_{n_k}/n_k^\alpha\to b $ a.s.\ as $k\to\oo$. 
It is interesting to investigate the relationship between $B(F,(n^\alpha))$ and 
its subset $\Bs(F,(n^\alpha))$ of \lq strong limits'.

\begin{theorem}\label{sl2}
{\bf\cite{MR359003}}
\begin{letlist}
\item Let $\alpha > \frac12$, $\alpha\ne 1$, and suppose $\Bs(F,(n^\alpha))\cap\RR\ne\es$.
Then $\Bs(F,(n^\alpha))\cap \RR$ equals one of
$(-\oo,0]$, $[0,\oo)$, $\RR$, and
each of these  three possibilities can occur.
\item Let $\alpha=1$. For any closed subset $B$ of $\ol\RR$, there exists
$F$ such that $\Bs(F,(n))=B$.
\end{letlist}
\end{theorem}

Harry  later verified a conjecture of Erickson \cite{erick} that a general $d$-dimensional random walk
diverges at least as fast as does a simple symmetric random walk. 

\begin{theorem} {\bf\cite{MR489585}}
Let $(X_i)$ be iid $\ZZ^d$-valued random variables whose support
is contained in no hyperplane, and let $S_n=\sum_{i=1}^n X_i$. 
Let $\psi:[1,\oo)\to(0,\oo)$ be such that $\psi(t)/t^{\frac12}\to 0$ as $t \to \oo$, and let
$(S_n^*)$ be a simple symmetric random walk on $\ZZ^d$.
If $|S_n^*|/\psi(n)\to \oo$ a.s.\ as $n\to\oo$, then
$|S_n|/\psi(n)\to \oo$ a.s.\ as $n \to\oo$.
\end{theorem}

In joint work with Durrett and Lawler \cite{MR1146451,MR1159577} published under the banner of gambling,
Harry investigated the long-term behaviour of a walk $(S_n)$ on $\ZZ$ with independent steps
$(X_i)$, each of which may 
have one of a given finite set $\sD$ of zero-mean distributions.  In an answer to a question
of Sp\u{a}taru, they showed that
any of the following may occur: 
(i) $\PP(S_n\to\oo)=1$, or (ii) $\PP(S_n\to-\oo)=1$, or
(iii) there are transient oscillations between $\pm\oo$.
Condition (i) is equivalent to $\PP(S_n>0)\to 1$, and similarly for (ii).
Now suppose $\sD=\{F_1,F_2\}$  where $F_i$ has infinite variance but finite $q_i$th moment, with $q_i\in(1,2)$.
When $q_1+q_2>3$, we have simultaneously that $\limsup_{n\to\oo}\PP(S_n<0)>0$
and $\limsup_{n\to\oo}\PP(S_n>0)>0$, whatever the choice of  the distributions of the $(X_i)$.

\subsection{Random walk after 1992}

Following his solution in 1980 to the percolation problem on $\ZZ^2$
(see Theorem \ref{pc12}), Harry's interests shifted towards percolation 
and related areas. He returned to more classical problems on 
a visit to Australia in 1990, and between 1992 and 2004 he published 10 papers with Ross Maller on random walks. 
This work is somewhat technical, and complete statements of their results (including
the modes of convergence) are not included here.

Let $(X_i)$ be an iid sequence with sums $S_n=\sum_{i=1}^n X_i$. 
For given $n \ge 1$, let $X_{(1)}\ge X_{(2)} \ge \dots\ge X_{(n)}$
denote the (decreasing) order statistics of the subsequence $(X_i:i=1,2,\dots,n)$, 
and let $X^{(1)}, X^{(2)}, \dots, X^{(n)}$
be the subsequence reordered in decreasing order of absolute values. The \lq trimmed sums' are defined as
$$
{}^{(r)}S_n = S_n-X_{(1)}-\dots - X_{(r)},
\qq
{}^{(r)}\ol S_n = S_n-X^{(1)}-\dots -X^{(r)},
$$
for $r=0,1,2,\dots,n$, and they express in two ways the degree of dominance of the
overall sum $S_n$ over the summands. Working with the trimmed sums for fixed $r$ is termed \lq light trimming'.
The results of
\cite{MR1188043} 
include necessary and sufficient conditions for the lightly trimmed sums
to satisfy $ {}^{(r)}S_n/X_{(r)} \to\oo$
and $ {}^{(r)}\ol S_n/|X^{(r)} |\to\oo$ as $n\to\oo$. In \cite{MR1455151}, 
Kesten and Maller studied the validity of the limits $S_n\to\oo$ and $S_n/n\to\oo$, and obtained
conditions under which $\oo\in \Bs(F,1)$ and/or $\oo\in \Bs(F,n)$, in the notation
of  Theorem \ref{sl2}. They turned in \cite{MR2057578} 
to the effect of trimming on a generalized law of the iterated logarithm.

In papers \cite{MR1303651,MR1337473}, Kesten and Maller established necessary and sufficient conditions for
 $ {}^{(r)}S_n/C_n \to\pm\oo$
and $ {}^{(r)}\ol S_n/C_n\to\pm\oo$ as $n\to\oo$, where the $C_n$ are non-random.
The results are of course related to Harry's earlier work on strong limits of 
the ratio $S_n/\gamma_n$; see 
Theorems \ref{sl1}--\ref{sl2}. They imply, for  example, that $\PP(S_n>0)\to 1$ if and only if
$S_n\to \oo$ in probability, a fact of interest to compulsive gamblers with deep pockets.

Moments of the first and last exit times of an interval $(-\oo,x]$
are explored in \cite{MR1408415}, and the orders of magnitude
of such moments are established as $x \to \oo$. The so-called \lq stability' of
exit times from a strip or half-plane of $\RR^2$ is studied in \cite{MR1725708},
with potential connections to sequential analysis.

In two linked articles \cite{MR1660924,MR1645198}, the authors investigated the ratio $S_n/n^\kappa$
as $n \to\oo$. Suppose that $\PP(S_n\to\oo)=1$, and  let $\kappa\ge 0$,
$\kappa\ne 1$, and when $\kappa>1$ assume also that $\EE(|X_1|^{1/\kappa})=\oo$. 
Then the probability that $S$ crosses the curve $y=ax^\kappa$ before it crosses $y=-ax^\kappa$ 
tends to $1$ as $a\to\oo$. This intuitively clear statement is not simple to prove, and 
indeed it is false when $\kappa=1$. These results are connected in the second paper to
the almost-sure limits $\limsup_{n\to\oo} S_n/n^\kappa=\oo$ and $\limsup_{n\to\oo} |S_n|/n^\kappa=\oo$,
for which necessary and sufficient conditions are given.

The \emph{range} $R_n$ of a random walk is the number of distinct points visited
up to time $n$. Hamana and Kesten \cite{MR1841327,MR1899229} 
proved a large-deviation principle for the range of a 
random walk on $\ZZ^d$ with $d \ge 1$. The subadditivity of the sequence $(R_n)$
is useful, and for once it turns out that the proof is hardest when $d=1$.  
Certain results and conjectures for the range of random walk on a Cayley graph are presented in \cite{MR2299929}.

\subsection{Random walk in random environment}\label{sec:rwre}

Random walks commonly live on either $\ZZ$ (or $\ZZ^d)$, or more generally on $\RR^d$ for some $d \ge 1$.
Groups provide a more general setting as in Section \ref{sec:rwgp}, and lately there has been great interest
in random walks on general connected graphs. In all such settings, the underlying graph 
$G$ is fixed,  the randomness is associated with the walk only, 
and the ensuing process constitutes a Markov chain.
The situation is much more complex if the environment is random also. 
In such a case, the position of the walker
is no longer a Markov chain, and there are generally sub-domains of the environment where
the walk moves anomalously. Such a process is termed `random
walk in random environment' (RWRE). 

On the one-dimensional line $\ZZ$, simple RWRE amounts to: 
(i) sampling random variables for each site, that prescribe the transition probabilities 
(left or right?) when the walker reaches the
site, and (ii) performing a random walk (or more precisely a
Markov chain) with those transition probabilities.  RWRE 
is a Markov chain \emph{given the environment} of transition probabilities, but it is not itself 
Markovian because, as the walker moves, it accrues information
about the random environment.

Let $A=\{\alpha_x: x=\dots,-1,0,1,\dots\}$ be iid random variables taking values in $[0,1]$,
and let $\beta_x=1-\alpha_x$. Conditional on $A$, 
let $X=(X_n)$ be a random walk on $\ZZ$ which, when at position $x$, moves one step rightwards (\resp,
leftwards) with probability $\alpha_x$ (\resp, $\beta_x$). 

Let $\rho=\beta_0/\alpha_0$. Suppose that there exists $\kappa\in(0,\oo)$ such that
$$
\EE\log \rho<0, \q \EE(\rho^\kappa)=1, \q 
\EE(\rho^\kappa\log^+ \rho)<\oo,
$$
and in addition $\log \rho$ is non-arithmetic. The following, amongst other
things,  was proved by Kesten, Kozlov, and Spitzer in 1975.

\begin{theorem} {\bf\cite{MR380998}}
There exist constants $A_\kappa,B_i>0$, and stable laws $L_\kappa$
with index $\kappa$, such that the following hold as $n\to\oo$, for all appropriate $x \in \RR$.
\begin{letlist}
\item If $\kappa<1$, $\PP(X_n \le xn^\kappa) \to 1-L_\kappa(x^{-1/\kappa})$.
\item If $\kappa=1$, for suitable $\delta_n \sim  n/(A_1\log n)$,
$$
\PP\bigl( X_n-\delta_n \le xn/(\log n)^2\bigr) \to1- L_1(-xA_1^2).
$$
\item If $1<\kappa< 2$, 
$$
\PP\bigl(X_n-(n/A_\kappa)\le xn^{1/\kappa}\bigr) \to 1- L_\kappa(-xA_\kappa^{1+1/\kappa}).
$$
\item If $\kappa=2$, 
$$
\PP\bigl(X_n -(n/A_2) \le x(n\log n)^{\frac12}B_1A_2^{-3/2}\bigr) \to \Phi(x).
$$
\item If $\kappa>2$,
$$
\PP\bigl(X_n -(n/B_3) \le xn^{\frac12}B_2B_3^{-3/2}\bigr) \to \Phi(x).
$$
\end{letlist}
Here, $L_\kappa$ is concentrated on $[0,\oo)$ if $\kappa <1$, and has mean zero when $\kappa>1$.
The standard normal distribution function is denoted $\Phi$.
\end{theorem}

In particular, when $\kappa<1$, $X_n$ has order $n^\kappa$.
The complexities of the main theorem of \cite{MR380998} are reproduced above in
illustration of Harry's determination and ability to get right to the heart of a problem.  

In the critical (recurrent) case, when  $\rho$ is bounded away from $0$ and $1$, and in addition satisfies
$$
\EE\log \rho=0, \qq \sigma^2 := \EE[(\log\rho)^2]\in(0,\oo),
$$
Sinai \cite{sinai} showed that $\sigma^2X_n/(\log n)^2 \Rightarrow L$
under the annealed measure $\PP$, for some functional $L$ of the standard Wiener process. Harry was able to calculate the density function of $L$ 
(see also \cite{Golosov}).

\begin{theorem} {\bf\cite{MR865247}}
The random variable $L$ has density function
$$
\frac d{dx}\PP(L\le x) = \frac2\pi\sum_{k=0}^\oo \frac{(-1)^k}{2k+1}
\exp\left\{-\frac{ (2k+1)^2\pi^2}8 |x|\right\},
\qq x\in \RR.
$$
\end{theorem}

Harry returned to this RWRE in \cite{MR0458648}, showing, subject to the condition
$\EE\rho<1$, that the environment seen from the position of the walker has a limit distribution.
Partial results when $\EE\rho=1$ were stated also, with proofs \lq\lq available from the author".

The above may be viewed as  early contributions to the theory
of RWRE, and to the related area of homogenization of differential operators with
random coefficients. 
RWRE in higher dimensions requires different techniques, 
and major progress has been made by Bricmont and Kupiainen \cite{brick} and others
(see also the references in Sznitman's paper \cite{sz}).

There are a number of variants of the RWRE problem.
One such variant arises when the walk lives on a random graph such as a tree or part of a lattice, and
Harry's work in these settings is summarised in Theorems \ref{thm:ao2} and \ref{subdiff2}--\ref{thm:rwinfc}.
In other systems, the environment evolves in time in a manner that depends on the trajectory of the random walk.
In the edge-reinforced random walk model, each edge is traversed with a
probability depending  on its weight, and its weight increases as it is traversed. With Durrett and Limic
\cite{MR1902191}, Harry showed that `once-reinforced' walk (ORRW) of \cite{davis}
on a tree is transient, and the distance of the 
walker from its starting point grows linearly with time. The ORRW model 
has attracted attention since, and the reader is referred to \cite{coll20} for recent results.

Finally, den Hollander, Kesten, and Sidoravicius \cite{MR3217035} studied a two-type particle system on $\ZZ^d$ in which
the jump-distribution  of a type-2 particle  depends on whether or not there is type-1 particle present.

Whereas the above is  set
in \emph{discrete} space--time, Harry's two papers \cite{MR528186,MR597029} with Papanicolaou 
from 1979/1980 are
directed at diffusion in a \emph{continuous} random environment.  In a model for turbulent diffusion, they considered
the stochastic differential equation $\dot x(t)=v+\epsilon  F(x(t))$ where  
$v\ne 0$ is a fixed vector and $F$ is a zero-mean stationary process satisfying a smoothness condition.
Their main result is that diffusivity holds in the limit as $\epsilon \downarrow 0$. Stochastic acceleration, and the equation
$\ddot x(t)=\epsilon F(x,\dot x)$ in $d\ge 3$ dimensions, is the topic of the second paper, where they 
conclude that  $x$ converges to a diffusion process as $\epsilon \downarrow 0$.

\subsection{Diffusion limited aggregation}\label{sec:dla}

\begin{figure}[h]
\centerline{
\includegraphics[width=0.5\textwidth]{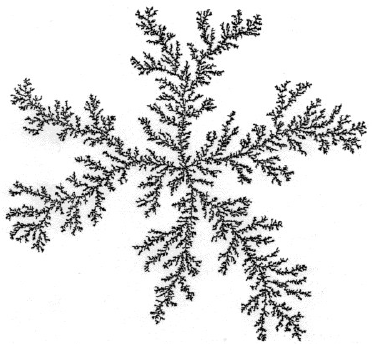}}
\caption{A 
simulation of a two-dimensional diffusion limited aggregation.
(By courtesy of Aliette Mandelbrot.)}
\label{fig:dla}
\end{figure}

Diffusion limited aggregation (DLA) is a stochastic model believed to exhibit dendritic growth
and power-law correlations, which was introduced in finite volumes by Witten and Sander \cite{wittens} (see Figure \ref{fig:dla}).   
We describe next Harry's results for a DLA model for a growing cluster on the lattice  $\ZZ^d$ with $d \ge 2$.

Let $A_0$ be the origin of $\ZZ^d$. The intuition of growth is as follows.
Conditional on the set $A_n$, the set $A_{n+1}$ is obtained by starting a
random walk \lq at infinity' and stopping it when it reaches
a point that is adjacent to $A_n$, and then adding that
point to the set $A_n$.  The concept of \lq release at infinity' 
may be interpreted as starting at a point  $x$ and then letting $|x|\to\oo$.
When $d=2$, this random walk hits $A_n$ with probability $1$, whereas when $d \ge 3$
the limit probability of hitting $A_n$ is $0$. In the latter case, one conditions on hitting 
$A_n$, and then passes to the limit as $|x|\to\oo$. This mechanism may be formalised using
so-called harmonic measure (\cite{MR915132,MR1146455}).

Computer simulations of \cite{wittens} and others 
indicate that $A_n$ grows in the manner of a fractal set,
with arms that extend across space. and whose embrace tends to trap the incoming particle when 
it is still far from  the origin.
One way of measuring this effect is via the radius of $A_n$, $\rad(A_n)=\max\{|a|: a\in A_n\}$,
and it is believed that $\rad(A_n)$ grows in the manner of $n^\alpha$ for some $\alpha$.
It is immediate by the isoperimetric inequality that such $\alpha$ satisfies $d^{-1} \le \alpha \le 1$.

\begin{theorem} {\bf\cite{MR873177,MR1077203}} \label{thm:dla}
There exists a constant $C=C(d)$ such that, a.s.\ eventually,
$$
\rad(A_n) \le \begin{cases} Cn^{2/(d+1)} &\text{if } d\ge 2,\ d \ne 3,\\
C(n \log n)^{1/2} &\text{if } d=3.
\end{cases}
$$
\end{theorem}

The $d=3$ bound has been improved slightly
by Lawler \cite[Thm 2.6.1]{lawler-int} to
$\rad(A_n) \le 
Cn^{1/2}( \log n)^{1/4}$. Subject to this, the results of
Theorem \ref{thm:dla} (dated 1990) remain the best rigorous results currently known
for DLA on $\ZZ^d$. No non-trivial lower bound on $\alpha$ is currently known. 

Simulations  suggest that $\rad(A_n) \sim C n^\alpha$ where the \lq dimension' $1/\alpha$ satisfies
$$
\frac1\alpha\approx \begin{cases}
1.7 &\text{if } d=2,\\
\dfrac{d^2+1}{d+1} &\text{when $d$ is large}.
\end{cases}
$$
See, for example, \cite{vicsek}.

\subsection{Random walk in random scenery}\label{sec:rwrs}

The subject of random walk in random scenery (RWRS) seems to have been initiated by 
Kesten and Spitzer \cite{MR550121}\footnote{Published in a volume of the 
\emph{Zeitschrift f\"ur Wahrscheinlichkeitstheorie und verw.\ Geb.},
and dedicated to its first editor,  Leopold Schmetterer, on his 60th birthday in 1979.}.
Let $\{\xi(x):x \in \ZZ\}$ be iid random variables, and let $W_n=\sum_{i=1}^n \xi(S_n)$
be the cumulative scenery of an independent random walk $(S_n)$. The 
process $(W_n)$, when suitably normalized, 
converges  to a self-similar process, under the assumption that the $\xi(x)$ and the 
steps $X_i$ belong to the domains of attraction of stable laws.  

A number of authors, including Kasteleyn \cite{kast85} and Keane and den Hollander \cite{KdH}, developed the
ergodic theory of the scenery process viewed progressively by the random walk.

In a more recently considered RWRS problem,
suppose there are only two possible sceneries $\xi$, $\xi'$
on a graph $G$. By progressive observation of its local scenery, 
can a random walker on $G=(V,E)$ decide which of the two is the true scenery? In a randomized version of this question,
the set of sceneries is taken as $\Xi(k)=\{0,1,\dots,k\}^V$, and 
a random scenery $\eta\in\Xi$ is picked according to uniform product measure on $\Xi(k)$. 
The concept of \lq distinguishability between sceneries' requires a definition, and once that is done,
one arrives at the following theorem of Benjamini and Kesten.

\begin{theorem}\mbox{\hfil}
\begin{letlist}
\item  {\bf\cite{MR1428097}}
For $G=\ZZ,\ZZ^2$ and $k=1$, and any fixed scenery $\xi\in\Xi(1)$, we have that 
almost every $\eta$ is distinguishable from $\xi$. 
\item {\bf\cite{MR1439524}} Let $G=\ZZ$, and let $\what \xi$ be obtained from the scenery $\xi\in \Xi(k)$ 
by altering the value $\xi(0)$ at the origin only. If $k \ge 5$, for almost every $\eta\in\Xi(k)$,
we have that $\eta$ and $\what \eta$ are distinguishable.
\end{letlist}
\end{theorem}

There has been a great deal of progress with this beautiful problem since the above,
due in part to Henry Matzinger, one of Harry's former PhD students. The reader
is directed to \cite{mpp} for relevant references.

\section{Products of random matrices}  \label{sec:prm}

Harry Kesten and Hillel Furstenberg spent the academic  year 1958--9 as instructors at Princeton University, 
where they became interested in recent work of Bellman. In 
the article \cite{bell}, one of 19 works by Bellman
listed on MathSciNet for 1954, he derived some results for the asymptotic behaviour of the product
of $n$ independent, random $2 \times2$ matrices. Furstenberg and Kesten were able to extend these results 
very considerably in their now well-known paper, which was described by Bellman as \lq\lq difficult and ingenious''.

\begin{theorem} {\bf\cite{MR121828}} \label{thm:kf}
Let $M_1,M_2,\dots$ be a stationary ergodic sequence of random $d\times d$ matrices,
such that $\EE \log^+\|M_1\|<\oo$. Then $Y_n=M_nM_{n-1}\cdots M_1$ satisfies
$$
\lim_{n\to\oo} \frac1n \log\|Y_n\|=E\qq\text{a.s.,}
$$
where
$$
E=\lim_{n\to\oo}\frac1n \EE\bigl(\log\|Y_n\|\bigr).
$$
\end{theorem}

Here, $\|A\|=\max_i\sum_j |a_{ij}|$ for $A=(a_{ij})$. A central limit theorem holds subject to
appropriate conditions on the $M_n$.

Their analysis made use of some subadditivity, but they did not anticipate
the forthcoming theory of subadditive stochastic processes (see Section \ref{sec:fpp}), 
which would later furnish a short proof of some of their results.
The Oseledec multiplicative ergodic theorem, which appeared a few years later (\cite{osel}),
may be viewed as a generalization of the results of \cite{MR121828}.

In the article \cite{MR440724}, Harry studied the
stochastic recurrence $Y_n=M_nY_{n-1}+Q_n$, where the $(M_n,Q_n)$ are 
independent, identically distributed pairs of $d\times d$ matrices, and the $Y_n$
are $d$-vectors. He noted that such equations occur in a variety of settings  including RWRE, branching processes in a random
environment, control theory, and models for evolution and cultural inheritance.
Key to this study is the behaviour of the products $M_1M_2\cdots M_n$. His results include
a connection to the stable laws, and a renewal theory for products of random matrices.
Near the heart of the work (when $d=1$)
lies the stochastic equation $Y\eqd MY+Q$, 
and the conclusion that such $Y$ can be heavy-tailed. 
This work has generated a considerable amount of interest since in probability, statistics, and mathematical finance, 
as noted by the authors of \cite{BDM}: \lq\lq The highly praised 
Kesten paper has motivated several generations of
researchers to work on closely related topics''.

Later, in work motivated by the potlatch process,
Kesten and Spitzer \cite{MR761563} investigated the weak convergence
of products of independent non-negative random matrices, and particularly conditions under  which
the product converges weakly, without normalization, to a measure that is not concentrated on the zero matrix.  

\section{Self-avoiding walks}\label{sec:saw}

A \emph{self-avoiding walk} (SAW) on a graph is a
path that visits no vertex more than once. SAWs were introduced by Orr \cite{orr} as a simple model
for long-chain polymers (see Flory \cite{flory}), and counting SAWs on a given graph is a fundamental
combinatorial problem with ramifications in both mathematics and physics.
 
Let $\ZZ^d$ be the $d$-dimensional hypercubic lattice with $d \ge 2$,
and let $\chi_n$ be the number of $n$-step SAWs starting at the origin.
The principal  counting problem is to establish the asymptotic behaviour 
of $\chi_n$ as $n\to\oo$, and to determine the 
radius of a \lq typical' SAW of length $n$. In a notable paper \cite{hammm}, Hammersley and Morton used
subadditivity to show that
that there exists a \lq connective constant' $\kappa=\kappa(\ZZ^d)$ given by 
$$
\frac 1n \log \chi_n\to\log\kappa\qq\text{as } n\to\oo,
$$ 
which is to say that $\chi_n=\kappa^{n(1+\o(1))}$. 

Self-avoiding walks have earned a reputation 
for being hard to study, and there remains much to understand despite progress over more than 50 years.
For example it remains an open problem to prove that
$$
\chi_{n+1}/\chi_n \to \kappa \qq\text{as } n\to\oo,
$$ 
despite Harry's ratio-limit theorem from 1963.

\begin{theorem}\label{thm:ksaw}
Let $d \ge 2$.
\begin{letlist}
\item {\bf \cite{MR152026}} The number $\chi_n$ of $n$-step SAWs on $\ZZ^d$ starting at the origin satisfies 
$$
\left|\frac{\chi_{n+2}}{\chi_n}-\kappa^2\right|\le An^{-\frac13},
$$
for some constant $A$.
\item {\bf\cite{MR166845}}
We have that $\chi_n \le \kappa^n\exp\{Bn^{2/(d+2)}\log n\}$ for some constant $B$.
\end{letlist}
\end{theorem}

Part (a) was proved using a surgery that replaces one part of a SAW by another; the $2$ appears
because, since  $\ZZ^d$ is bipartite, parity of lengths is preserved in surgery. In the proof,
Harry introduced an argument known now as \lq Kesten's pattern theorem', which has since been very useful
in a number of contexts. It states that
any configuration of $k$ consecutive steps, that has the property that it can 
occur more than once in an $n$-step SAW, has to occur at least $an$ times, for some $a>0$, 
in all but \lq exponentially few' such SAWs. 

While it is a simple consequence of subadditivity that $\chi_ n \ge \kappa^n$, upper bounds on $\chi_n$ are harder to derive.
Hammersley and Welsh \cite{hammw}
had earlier proved that $\chi_n\le \kappa^n \exp\{Bn^{\frac12}\}$ for all $d \ge 2$, and
Theorem \ref{thm:ksaw}(b) was an improvement of this for $d \ge 3$. The Hammersley--Welsh inequality remained the
best available when $d=2$ until the appearance in 2018 of \cite{dghm} with an improved exponent
of order $n^{\frac12-\epsilon}$.
In 1992, Hara and Slade \cite{haras2,haras1} developed a lace expansion for SAWs that yielded
$\chi_n \sim C\kappa^n$ when $d \ge 5$. Part (b), above,  remains the best current  upper bound
when $d=3,4$. As remarked in \cite[p.\ 69]{mads}, Harry 
stated part (b), but he proved only a weaker bound  in \cite{MR166845}
since he hoped someone would find a sharper inequality.

The theory of self-avoiding walks on two-dimensional lattices is an important topic
of current research.
Duminil-Copin and Smirnov \cite{dumcs} proved the longstanding conjecture that the connective constant
of the hexangular lattice $\HH$ satisfies $\kappa(\HH)=\sqrt{2+\sqrt 2}$, and in so doing
presented a rigorous connection between two-dimensional
SAWs and conformal invariance. It is a notorious open problem
to prove that, when suitably rescaled, a uniformly random $n$-step SAW on $\HH$ converges 
as $n\to\oo$ to so-called SLE$_{8/3}$.

\section{Branching processes}\label{sec:bp}

The branching process (or, as Harry liked in later years to write, the Bienaym\'e--Galton--Watson process)
is the most fundamental stochastic model for population growth. Each individual has a family
of offspring, and the family-sizes are drawn independently from a given distribution. Branching processes
have been used in a multitude of application areas, and they continue to be a source of wonderful mathematical problems.  

Let $Z_n$ be the size of the $n$th generation of a branching process, and suppose $Z_0=1$
and $\mu=\EE(Z_1)<\oo$. It is
easy to see that $\EE(Z_n)=\mu^n$, so that $\mu$ may be viewed as the average rate of population growth.
It is an undergraduate exercise (due to Christensen and Steffensen in 1930, see \cite{albertsen}) 
that the process $Z$ becomes a.s.\  extinct, in that a.s.\ $Z_n=0$ for large $n$, 
if $\mu<1$, and not if $\mu>1$. In the latter case, the extinction probability $\eta$ is
easily calculated. This and very much more is explained in Harris's influential book
\cite{Har-BP} from 1963 on the subject, which includes an account of \lq multitype' processes with more than one type of individual. 

Harry wrote a  number of significant articles about branching processes.
Perhaps his widest known work is his paper \cite{MR198552}, one of three joint works with Stigum
concerned with multitype branching processes, written around 1966,
namely \cite{MR200979,MR198552,MR205340}.  They extended work of Levinson
\cite{levin} concerning the convergence of supercritical multitype branching processes, and gave
a necessary and sufficient condition concerning the mean of the associated martingale $W_n=Z_n/\mu^n$.
This condition was new even for processes with only one type, and it is in this 
simpler form that it is presented here.

\begin{theorem} {\bf\cite{MR198552}}
Let $Z=(Z_n)$ be a branching process with $Z_0=1$ and  $\mu:=\EE(Z_1)>1$. 
The almost-sure limit $W=\lim_{n\to\oo} W_n$  satisfies $\EE(W)=1$ if and only if
$\EE(Z_1\log^+ Z_1)<\oo$. Under this condition,  $W$ has an atom of size $\eta$ at $0$,
and is absolutely continuous on $(0,\oo)$.
\end{theorem}

These conclusions were known already subject to a second moment condition, and the authors of \cite{MR198552} 
observed that the main novelty of their work
lay in the identification of 
the condition $\EE(Z_1\log^+Z_1)<\oo$. An interesting discussion of these and further results is found in \cite{lpp}.

Celebrated theorems of Kolmogorov \cite{kolm1} and Yaglom \cite{yag} identified the speed of extinction of
a critical branching process (when $\mu=1$), and its weak limit conditional on non-extinction.
With Ney and Spitzer, Harry weakened their third-moment condition
to the necessary and sufficient second-moment condition, using \lq\lq traditional" tools.
 
\begin{theorem}[\cite{MR0207052}]
Let $Z=(Z_n)$ be a branching process with $Z_0=1$ and  $\mu:=\EE(Z_1)=1$.  
If $\sigma^2:= \var(Z_1)$ satisfies $0<\sigma^2 \le \oo$, then
$$
n\PP(Z_n>0) \to \frac 2{\sigma^2},\qq
\PP(Z_n > nx \mid Z_n>0) \to e^{-2x/\sigma^2},
$$
as $n\to\oo$.
\end{theorem}

The incipient infinite percolation cluster $\wt C$ is a connected subgraph of $\ZZ^2$ with fractal-like qualities.
One way to understand something about  the geometry of $\wt C$ is to study the variance of 
random walk thereon,
and in particular to show that it is subdiffusive (see Section \ref{sec:incic}). Perhaps
as a warm-up to this problem, Harry included in his paper \cite{MR871905}
an account of random walk on the family tree of a critical branching process conditioned on non-extinction,
and his results include the following.

\begin{theorem} {\bf\cite{MR871905}} \label{thm:ao2}
Let $Z$ be a critical branching process satisfying $Z_0=1$ and $0<\var(Z_1)<\oo$,
and let $T$ be a rooted tree with the law of $Z$ conditioned on non-extinction.
Let $X=(X_n)$ be a random walk on $T$, and write $\PP$ for the annealed law of $X$.
The height $h(X_n)$ is such that $n^{-\frac13}h(X_n)$
converges weakly (under $\PP$) to a distribution with no mass at $0$.
\end{theorem}

The quenched version of the above statement is left open.
A corresponding  invariance principle is stated but not proved in \cite{MR871905}, since
\lq\lq at this time we only have a monstrously long proof of [$\cdots$] and we
therefore restrict ourselves to the following weaker result $\dots$".
The $n^{\frac13}$ scaling is consistent with the so-called Alexander--Orbach conjecture
in the setting of a critical branching process conditioned on non-extinction (\cite{AlexO}).  

Harry continued his work on the family-tree $T$ of a critical branching process in several papers,
including \cite{MR494543} on
branching Brownian motion,  and
\cite{MR1311717} on an application to river basin hydrology
published in the Dynkin festschrift.
Kesten and Pittel \cite{MR1603252} established an asymptotic formula for the joint
distribution for the height and number of leaves of $T$. 

With Durrett and Waymire in \cite{MR1088403},
Harry studied the family-tree $T_n$ of a branching process, when conditioned to have total size $n$. 
The issue in question was to understand the growth 
of the weighted height of $T_n$, when iid weights are associated with the edges of $T_n$.
They explored the manner in which the answer depends on the power-law of the weight distribution.
Harry returned in \cite{MR1353560} to this edge-weighted process, viewed as a critical branching random walk
in which a vertex $w$ is displaced from its parent $v$ by the weight of the corresponding edge
$\langle v,w\rangle$. Let $M_n$ be the maximal displacement of the $n$th generation. 
If the edge-weights have mean $0$ then, conditional on non-extinction at time $\beta n$
where $\beta\in(0,\oo)$,  the sequence $n^{-\frac12}M_n$
converges in distribution. 
This is related to results for Aldous's continuum random tree and Le Gall's Brownian snake.

Harry considered  in \cite{MR343386} a probabilistic model for \lq\lq zygotic 
mutation without selection at one locus'' in population genetics. Subject to a quite technical assumption, he proved a convergence theorem 
for this model which generalises the principal limit theorem for a supercritical branching process.

Estimates concerning the geometry of a branching process turn out to be key to Harry's study with Grimmett
of the electrical resistance of a complete graph (\cite{MR760886,0107068}).
In \cite{MR1031281}, he developed detailed asymptotics for a supercritical branching process 
with a countable infinity of types, in order to derive properties of projections of random Cantor sets.

\section{Percolation theory}\label{sec:perc}

\subsection{The critical probability in two dimensions}

\begin{figure}
\centerline{\includegraphics[width=0.65\textwidth]{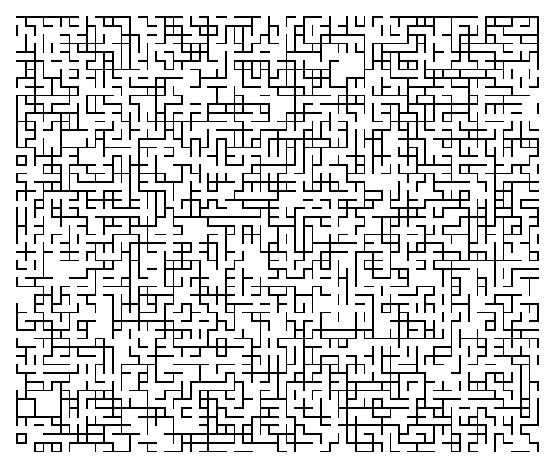}}
\caption{A simulation of bond percolation on $\ZZ^2$ with $p=0.51$.}
\label{fig:perc}
\end{figure}

The percolation process is a canonical model for a disordered medium.
For concreteness, consider the $d$-dimensional hypercubic lattice $\ZZ^d$ with $d \ge 2$,
and let $p \in [0,1]$.
We declare an edge of $\ZZ^d$ to be \emph{open} with probability $p$, and \emph{closed} otherwise,
different edges having independent states. 
What can be said about the geometry of the subgraph induced by the open edges
(as illustrated in Figure \ref{fig:perc}), and
in particular for what values of $p$ does this graph 
possess an infinite cluster?

As Harry reflected in the preface of \cite{MR692943}, percolation is both a source of fascinating problems, 
and a central topic in the physical theory of phase transitions. 
The mathematical theory originated in work of Broadbent and Hammersley
\cite{BH57} on diffusion through a random medium, as a model for the transmission of particles through a face mask.
Percolation theory is now a key tool in the mathematics of interacting systems, and a
testbed for new techniques. The two-dimensional theory has been
especially prominent in the currently active theory of random planar geometry.

The \emph{percolation probability} is the function
$$
\theta(p)=\PP_p\bigl(\text{the origin lies in an infinite open cluster}\bigr),
$$
and the \emph{critical probability} is given by
$$
\pc(\ZZ^d)=\sup\{p: \theta(p)=0\}.
$$
It is fundamental that $0<\pc(\ZZ^d)<1$ when $d \ge 2$ (see \cite{Gr99} for a general account of percolation).
The obvious problem of the exact calculation of $\pc(\ZZ^d)$ is intractable when $d \ge 3$. 
The case $d=2$ is, however, very special, because of
a property of self-duality possessed by the square lattice, which motivated the 
then notorious conjecture that $\pc(\ZZ^2)=\frac12$.  Harry  proved this in 1980.

\begin{theorem} {\bf{\cite{MR575895}}} \label{pc12}
The critical probability of bond percolation on the square lattice is $\pc(\ZZ^2)=\frac12$.
\end{theorem}

The conjecture that $\pc(\ZZ^2)=\frac12$ went back at least to 1960. In 
that year, Harris \cite{teh} proved that $\theta(\frac12)=0$, whence $\pc(\ZZ^2)\ge \frac12$, using a 
clever geometrical
argument combined with a  positive-correlation inequality, both of which have lasting significance.
Hammersley's numerical simulations indicated a value strictly less than $\frac12$ (\lq\lq what better evidence could there
be for $\pc(\ZZ^2)=\frac12$'', he would say). There were important papers by Russo and Seymour/Welsh
that led to Harry's famous article \cite{MR575895} where he quantified the idea of \lq sharp threshold'
and hence derived a proof of Theorem \ref{pc12}.    
The related concept of \lq influence' has since been independently 
systematised in general settings by Kahn, Kalai, Linial \cite{KKL}, Talagrand \cite{T94}, and their 
successors (see \cite[Chap.\ 4]{g-pgs}).

Harry's subsequent volume \cite{MR692943} is a fairly formidable work on percolation in two
dimensions, designed unapologetically for the mathematician. 
These two works \cite{MR575895,MR692943} brought prominence to
percolation, and ushered in a period of intense activity, initially in two dimensions and later more generally.
The exact calculation of Theorem \ref{pc12} has been vastly extended since to percolation
on a large family of isoradial graphs
(\cite{gm2014}), that is, graphs embedded in $\RR^2$ in such a way that every face is convex with 
a circumcircle of given radius.
Such exact values were later extended to the critical points of random-cluster models on isoradial graphs (\cite{dclm}).

\subsection{Disjoint occurrence and the van den Berg--Kesten (BK) inequality}\label{sec:bk}

Correlation inequalities play a very important role in the theory of disordered systems in mathematics and physics.
The most prominent such inequality for percolation is the so-called (Harris--)FKG inequality, which states that
increasing events are positively associated. In investigating a converse inequality, van den Berg made 
a conjecture which he proved jointly with Harry.

Let $|E|<\oo$ and $\Omega=\{0,1\}^E$, endowed with a product measure $\PP$.
Let $A$, $B$ be increasing events in the finite product space $(\Omega, \PP)$. 
For $\om\in\Omega$ and $F \subseteq E$,
let $\om_F$ denote the configuration which agrees with $\om$ on $F$ and equals $0$ elsewhere.
The \emph{disjoint occurrence} event  $A \circ B$ is defined by
$$
A\circ B= \Bigl\{\om\in \{0,1\}^E: \exists\text{ $F\subseteq E$ such that 
$\om_F\in A$, $\om_{E\setminus F}\in B$}\Bigr\}.
$$

\begin{theorem}\label{thm:BK} {\bf\cite{MR799280}}
For increasing events $A, B\subseteq \Omega$, we have that $\PP(A \circ B) \le \PP(A)\PP(B)$.
\end{theorem}

This \lq BK inequality' is a delicate and tantalising result, with a number of useful applications.
It was extended by Reimer \cite{reimer} to a disjoint-occurrence inequality without the assumption that
the events be increasing. Van den Berg and Jonasson \cite{vdBJ} have proved Theorem \ref{thm:BK}
with $\PP$ replaced by the uniform measure on subsets of $E$ with given fixed cardinality.   

\subsection{Uniqueness of the infinite cluster}\label{sec:uic}

When $p>\pc$, there exists a.s.\ at least one infinite open cluster, but how many? This uniqueness problem
was answered by Harry in  joint  work with  Aizenman and  Newman. 
The proof is very general, but we state it here for the bond model on $\ZZ^d$.

\begin{theorem} {\bf\cite{MR901151}}
Consider bond percolation on $\ZZ^d$ with $d \ge 2$, and let $p\in[0,1]$.
The number $N$ of infinite open clusters satisfies
$$
\text{either} \q \PP_p(N=0)=1 \q \text{or} \q \PP_p(N=1)=1.
$$
\end{theorem}

Paper \cite{MR901151}
was soon superseded by the elegant argument of  Burton and  Keane \cite{bk1989}, but it remains important
as a source of quantitative estimates, as, for example, 
for the lower bound  of \cite{cerf2015} on the so-called two-arm exponent.  

Harry worked also on uniqueness in
long-range models on subgraphs of $\ZZ^d$. When $d=1$,
he and Durrett  \cite{MR1117011} extended a special calculation of Shepp \cite{shepp} in which sites $i$ and $j$ 
are connected with probability $p(i,j)$, where $p$ is a homogeneous function of order $-1$.
Turning to (spatially homogeneous) long-range percolation, let $p:\ZZ^d\to[0,1]$ be symmetric.
An undirected edge is positioned between sites $i$ and $j$ with probability $p(j-i)$. When is the 
ensuing random graph a.s.\ connected? The necessary and sufficient condition 
according to Grimmett, Keane, and Marstrand \cite{GKM} is 
that $\sum_z p(z)=\oo$ (subject to a natural irreducibility assumption on 
the support of $p$). The argument used there fails for subspaces of the form 
$\ZZ^{d-e}\times \ZZ_+^e$ with $e\ge 1$, and likewise the methods of \cite{MR901151,bk1989,KW}
do not seem to apply  for $d\ge 2$, $e\ge 1$.  
Harry showed in \cite{MR1188713}
that a corresponding summability condition is both necessary and sufficient.

\subsection{Incipience and random walk}\label{sec:incic}

Whereas Harris proved in \cite{teh} that $\theta(\frac12)=0$ when $d=2$, the proof 
that $\theta(\pc(\ZZ^d))=0$ for all $d \ge 3$ has 
so far resisted allcomers (it is now known for $d \ge 11$, \cite{FHofstad}, using 
an elaboration of the lace expansion of \cite{HaraS}).
There was interest among physicists in the idea that a large two-dimensional cluster is \lq nearly infinite'
when $p=\frac12$. Harry's response to this speculation was a pair of coordinated limit theorems.

\begin{theorem} {\bf\cite{MR859839}}
Consider bond percolation on $\ZZ^2$, and let $C$ denote the open cluster containing the origin $0$.
The following two conditional probability measures converge, and they have the same limit $\nu$:
\begin{numlist}
\item $\PP_{\frac12}(\cdot \mid K_n)$ as $n \to\oo$, where $K_n$ is the event that $0$ is connected 
by an open path to some lattice point with distance $n$ from $0$,
\item $\PP_p(\cdot\mid \text{\rm $C$ is infinite})$ as $p \downarrow \frac12$.
\end{numlist}
\end{theorem}

The measure $\nu$ (and the open cluster $\wt C$ at the origin under $\nu$) is called
the \emph{incipient infinite cluster}. While the fractal nature of $\wt C$ is still not properly understood, Harry
exhibited its mass dimension in terms of that of the $\PP$-probability $\pi_n$ that the origin is
connected to the half-space $[n,\oo)\times \ZZ$.

\begin{theorem} {\bf\cite{MR859839}} \label{massdim}
As $\epsilon\downarrow 0$, uniformly in $n$,
$$
\nu\left(\epsilon< \frac{|\wt C \cap B_n|}{n^2\pi_n} < \frac1\epsilon \right)\to 1,
$$
where $B_n=[-n,n]^2$.
\end{theorem}

Theorem \ref{massdim}  is complemented in \cite{MR859839} by the inequalities
$$
C_1n^{-\frac12+\eta_1} \le \pi_n \le C_2n^{-\eta_2},
$$
for some $\eta_i>0$. The sizes of large open clusters of $B_n$ were studied in \cite{MR1868996} by
Borgs, Chayes, Kesten, and Spencer for an interval of values of $p$ around the critical value  $\pc$,
with similar results for $3\le d\le 6$ subject to an extensive set of hypotheses. 

De Gennes proposed in \cite{deG} the metaphor of a random ant in order to explore
the geometry of a random set. In  a continuation of the above work, 
Harry studied a random walk $X=(X_n: n \ge 0)$ on the
incipient infinite cluster $\wt C$ starting at the origin.

\begin{theorem} {\bf\cite{MR871905}} \label{subdiff2}
For some $\epsilon>0$,  $(n^{-\frac12+\epsilon}X_n:n\ge 0)$ is tight. 
\end{theorem}

In particular, there exists $\epsilon>0$ such that
$n^{-\frac12+\epsilon}X_n \to 0$ as $n\to\oo$.
Since a diffusion process grows in the manner of $n^{\frac12}$, this implies the
subdiffusivity of $X$. This slowing occurs because of the dangling regions 
and blind allies on a multiplicity of scales
that are found in the random graph $\wt C$. More detailed results are available in the related situation 
of a random walk on the family tree of a critical branching process conditioned on non-extinction;
see \cite{MR871905} and Theorem \ref{thm:ao2}.

The incipient infinite cluster, with the  random walk thereon, has attracted much attention in a variety
of further contexts since Harry's work on $\ZZ^2$.
See, for example, the book of Heydenreich and van der Hofstad \cite{HeyvdH}
and the references therein. A well-known 
prediction of Alexander and Orbach \cite{AlexO}
states that, for random walk $X$ on the incipient  infinite cluster of $\ZZ^d$ (if it exists), we have
$$
\PP(X_{2n}=0) \approx n^{-\frac23}\qq \text{as } n\to\oo.
$$
This is believed to hold when $d \ge 6$, but has so far been proved 
(by Kozma and Nachmias) only in the large-$d$ case, which these
days means $d \ge 11$ (see \cite{NachK} and \cite{FHofstad,HeyvdH}).

De Genne's \lq ant in a labyrinth' is relevant also to random walk 
on the infinite open (supercritical) percolation cluster, 
where it provided a concrete example of a non-elliptic RWRE.
With Grimmett and Zhang, Harry showed the following early result in the area.

\begin{theorem} {\bf\cite{MR1222363}}\label{thm:rwinfc}
Let $d \ge 2$ and $p>\pc(\ZZ^d)$. Random walk on the a.s.\ unique infinite open cluster
is recurrent if and only if $d=2$.
\end{theorem}
  
Precise estimates have been obtained since for the transition probabilities and central limit theory
of this random walk, and the reader is referred to the papers \cite{ABDH,mb04} for accounts of 
recent work.

\subsection{Supercritical percolation}\label{sec:super}

The size $|C|$ of the cluster at the origin has an exponentially decaying tail in the subcritical phase when  $p<\pc$.
The situation is different when $p>\pc$, in that the distribution of  $|C|$ is in part controlled by the
size of its surface. The isoperimetric inequality suggests a
stretched exponential distribution.
There was a  proof by  Aizenman,  Delyon, and  Souillard \cite{ADS}
that $\PP_p(|C|=n) \ge e^{-\beta n^{(d-1)/d}}$ for some $\beta(p)>0$.
Kesten and Zhang \cite{MR1055419} showed the complementary inequality,
by a novel block argument that has been useful since.

\begin{theorem} {\bf\cite{ADS,MR1055419}} \label{ads}
Let $d \ge 3$ and $\pc(\ZZ^d)<p<1$. There exist functions $\beta$, $\gamma$,
taking values in $(0,\oo)$, such that
$$
e^{-\beta n^{(d-1)/d}} \le \PP_p(|C|=n) \le   e^{-\gamma n^{(d-1)/d}}, \q n \ge 1.
$$
\end{theorem}

When $d\ge 3$, the authors of \cite{MR1055419} were in fact only 
able to show the above upper bound for $p$ exceeding a certain value $\pslab$, but 
the full conclusion for $p>\pc$ followed once  Grimmett and  Marstrand \cite{GM} had proved 
the slab limit  $\pc=\pslab$.
Sharp asymptotics and the associated Wulff contruction were established later by  Alexander,  Chayes, and  Chayes
\cite{ACC} when $d=2$, and by  Cerf \cite{cerf00,cerf-stf}
in the substantially more challenging  situation of $d = 3$. 
The Kesten/Zhang construction of \cite{MR1055419} has recently proved key to the proof
of analyticity for the percolation probability $\theta$ in $d \ge 3$ dimensions; see \cite{gp20}. 

\subsection{Scaling theory}\label{sec:scal}

The critical probability $\pc$ marks a point of singularity separating the \emph{subcritical phase} 
(all open clusters are finite) from the \emph{supercritical
phase} (there exists an infinite open cluster). The singularity at $\pc$ has much in common with the
phase transitions of statistical physics,

Let $C$ denote the open cluster containing the origin.
According to scaling theory, macroscopic functions such as the percolation probability 
$\theta(p)=\PP_p(|C|=\oo)$ and the mean cluster size
$\chi(p)=\EE_p|C|$ have singularities at $\pc$ of the form $|p-\pc|$ 
raised to an appropriate power called a \emph{critical exponent}.
More generally, it is believed that
\begin{alignat*}{4}
    \theta(p)&\approx (p-\pc)^\beta,          &&\q \chi(p)&&\approx (p-\pc)^{-\gamma},\qq&&\text{as } p \downarrow \pc,\\
\kappa'''(p) &\approx |p-\pc|^{-1-\alpha}, &&\q \xi(p) &&\approx |p-\pc|^{-\nu},   \qq&&\text{as } p \to \pc,\\
&&&&&\llap{$\dfrac{\EE_p(|C|^{k+1}; |C|<\oo)}{\EE_p(|C|^k;|C|<\oo)}$}\approx |p-\pc|^{-\Delta},
   \qq&&\text{for $k \ge 1$, as } p \to \pc,
\end{alignat*}
where $\kappa(p)$ is the mean number of clusters per site, and $\xi(p)$ is the correlation length.
The asymptotic relation $\approx$ is to be interpeted logarithmically
(stronger asymptotics are also expected to hold).

\begin{figure}
\centerline{\includegraphics[width=0.8\textwidth]{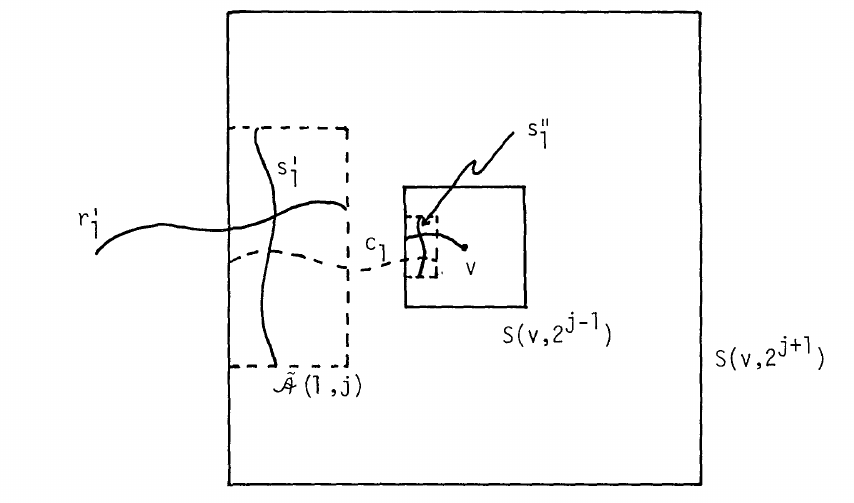}}
\caption{An example of a percolation figure by Harry Kesten, taken from his paper \cite{MR879034} on 
scaling theory.}
\label{fig:scaling}
\end{figure}

In similar fashion, when $p=\pc$, several random variables 
associated with the open cluster at the origin are expected to have power-law tails.  More specifically,
\begin{gather*}
\PP_{\pc}(|C|=n)\approx n^{-1-1/\delta},\q
\PP_{\pc}(\rad(C)=n)\approx n^{-1-1/\rho},\\
\PP_{\pc}(0\lra x)\approx |x|^{2-d-\eta},
\end{gather*}
as $n, |x| \to\oo$.

The set of critical exponents (including the eight given above) describes 
the nature of the singularity, and they are supposedly \lq universal' in that their values
depend on the number $d$ of dimensions but not on the choice of lattice. 
They are expected in addition to satisfy the so-called \lq scaling relations' of statistical physics. It is an open problem
to prove almost any of the above in general dimensions. 

Extraordinary (but incomplete) progress has been made towards the above when $d=2$. Harry's work has
been a key part of this program. Article \cite{MR633715} describes early progress
in the rigorous study of (non-)smoothness of macroscopic functions. Probably his most important
paper in this context was \cite{MR879034} on the scaling relations.

\begin{theorem} {\bf\cite{MR894549,MR879034}} \label{thm:scaling}
Let $d=2$. If the limits defining $\delta$ and $\nu$ exist, then
$$
\beta=\frac{2\nu}{\delta+1},\q  \gamma=2\nu\frac{\delta-1}{\delta+1},\q
\Delta=2\nu\frac\delta{\delta+1}, \q \eta=\frac 4{\delta+1}, \q 2\rho=\delta+1.
$$
\end{theorem}

This list of relations may seem somewhat sterile when stated in this bald manner.
Underlying these relations is a method for relating near-critical and critical percolation, supported
by a number of complex geometrical estimates of probabilities (somewhat systematised in
\cite{nolin}). These methods,
together with Theorem \ref {thm:scaling} itself, are an important part of the proof  
of the exact values of critical exponents for site percolation on the triangular lattice $\TT$:
$$
\beta=\tfrac5{36},\q \delta=\tfrac{91}5,\q\gamma=\tfrac{43}{18}, \q 
\eta=\tfrac5{24}, \q \nu=\tfrac43.
$$
This, and much more, has been proved over the last 20 years by 
Smirnov, Schramm, Lawler, Werner and others (see, for example, \cite{ww2009}).
The proof of the existence of critical exponents (for site percolation on $\TT$)
had to wait for the invention of the SLE processes by Schramm \cite{odeds};
Smirnov proved Cardy's formula, and explained conformal invariance (see \cite{ss-cras,cam-new});
and, finally, the team of Lawler, Schramm, and Werner built a theory of SLE and certain related lattice systems
(see, for example, \cite{os-books}).

The impact of Harry's papers  \cite{MR894549,MR879034} has been enormous across
two-dimensional percolation and related \lq percolative' systems including models such as dynamical,
frozen, and  invasion percolation'; see, for example, \cite{GPS,BKN,DSV}.

\subsection{First-passage percolation}\label{sec:fpp}

In this time-dependent percolation model on $\ZZ^d$,  the edges are assigned independent non-negative 
\lq time coordinates', interpreted as the time required to pass along the edge.
One studies the region reached from the origin within a given time.
This process was introduced by Hammersley and Welsh \cite{hw}
in 1965. They proved that, subject to certain conditions, the \emph{first passage time} $a_{0,n}$ from the origin to the point $(n,0,\dots,0)$
satisfies $a_{0,n}/n \to \mu$ for some \lq time constant' $\mu$. In so doing,
they introduced the notion of stochastic subadditivity, and they proved a version of the subadditive ergodic theorem.
Subadditivity became recognised as an important technique in probability and combinatorics, 
and it was a perennial theme in Harry's work. Indeed, it provides an elegant proof of 
part of his result with Furstenberg, Theorem \ref{thm:kf},
on the products of random matrices.
 
Harry and his co-authors resolved a number of significant problems in first-passage percolation,
and posed others, in a series of papers spanning nearly 20 years.  He established several fundamental properties
of the time constant in two dimensions, including positivity and continuity. 
Let $U$ be the law of the time coordinates, 
and $\mu=\mu(U)$ the time constant.

\begin{theorem}
Consider first-passage percolation on $\ZZ^2$.
\begin{letlist}
\item {\bf \cite{MR588406}}
If $U(0)<\frac12$, then $\mu(U)>0$.
\item {\bf\cite{MR633228}}
If $U_k \Rightarrow U$ in the sense of weak convergence, then $\mu(U_k)\to\mu(U)$.
\item {\bf\cite{MR1202515}} Let $U_1$ and $U_2$ be probability measures on $[0,\oo)$
with finite means such that $U_2$ is \lq\lq useful''. If $U_1\lest U_2$ but $U_1\ne U_2$, then $\mu(U_1)<\mu(U_2)$.
\end{letlist}
\end{theorem}

Concerning  part (a), it is immediate that $\mu(U)=0$ if $U(0)>\frac12$. Matters are more delicate when 
$U(0) = \frac12$, but nevertheless $\mu(U)=0$, and moreover Kesten and Zhang
\cite{MR1431216}  
proved a central limit theorem for the sequence $(a_{0,n})$, when suitably normalised, using 
martingale central limit theory. This work has been continued by several authors, including in the recent papers \cite{auffdh,damron2015,damron2017}.

Part (b), which was joint with Cox, was extended only in 2017 by Garet, Marchand,  Procaccia, and Th\'eret \cite{gmp} to $\ZZ^d$ with $d \ge 3$.
For an explanation of the \lq usefulness' condition of part (c), the reader is referred to the original paper. 
It implies in particular that $\mu(U_2)>0$. Stochastic ordering is denoted by $\lest$.

Mention should also be made of \cite{MR1221154}, which addresses the problem of the 
rate of convergence of $a_{0,n}/n$,
and of the degree of roughness of the boundary of the region reached from the origin up to time $n$. 
The scaling theory of such roughness remains an important open problem which has been considered
since by several authors (see \cite{chatt,dam18,MR3014795}). 
It is generally a hard problem to obtain decent lower bounds for a subadditive stochastic process
(see the relevant comments in \cite[p.\ 140]{auffdh}).

Exponential estimates and box arguments were a regular theme of Harry's work, and they featured
in his large-deviation theorem for first passage times.

\begin{theorem} {\bf\cite{MR751574}}\label{thm:ld}
Consider first-passage percolation on $\ZZ^d$ with $d \ge 2$, and
let $\epsilon>0$.
\begin{letlist}
\item If $U$ has finite variance, there exist $A,B>0$ such that
$$
\PP\bigl(a_{0n}<n(\mu-\epsilon)\bigr) \le Ae^{-Bn}, \qq n \ge 1.
$$
\item If $\int e^{\gamma x}\,U(dx)<\oo$ for some $\gamma>0$,  there exist $C,D>0$ such that
$$
\PP\bigl(a_{0n}>n(\mu+\epsilon)\bigr) \le Ce^{-Dn}, \qq n \ge 1.
$$
\end{letlist}
\end{theorem}

More refined analysis shows that the upper large deviations of part (b) can decay exponentially
in $n^d$, whereas the order of the lower deviations given above is typically correct. 
The contrast between lower and upper large deviations for subadditive functionals has attracted
quite a lot of attention since Harry's work. See, for example, \cite{chow-z,CGM}. 

Theorem \ref{thm:ld} has applications to the theory of random electrical networks, and to the diffusion constant in
a random medium. This connection led Harry to a fairly systematic study of 
the \lq random flow' problem: what can be said about the maximal flow across
a cube of $\ZZ^d$ with random capacities assigned to the edges.
For the sake of definiteness, consider the box $A_{n,k}=[0,n]\times[0,k]^2$ in $\ZZ^3$,
with bottom $B_{n,k}=\{0\}\times[0,k]^2$ and top $T_{n,k}=\{n\}\times[0,k]^2$.
To the edges $e$ of $A_{n,k}$ are assigned iid capacities $C(e)\ge 0$.
Let $\Phi_{n,k}$ be the maximum flow from $B_{n,k}$ to $T_{n,k}$ inside
$A_{n,k}$ subject to the constraint that the absolute value of the flow along any $e$ is no larger than $C(e)$.
The following is an exemplar of Harry's results.

\begin{theorem} {\bf\cite{MR876084,MR869483}}
Let $n=n(k) \to\oo$ as $k\to\oo$ in such a way that, for some $\delta>0$,
$k^{-1+\delta}\log n\to 0$. There exists $p_0>0$ such that, if $\PP(C=0)< p_0$ 
and $\EE(e^{\gamma C})<\oo$ for some $\gamma>0$, then the limit
$$
\nu=\lim_{k\to\oo} \frac 1{k^2} \Phi_{n,k}
$$
exists a.s.\ and in $L^1$.  Furthermore, $\nu\in(0,\oo)$.
\end{theorem}

Some tricky geometrical issues arise in the proof,
concerning  the combinatorics and topology
of dual surfaces in $\RR^3$. The condition involving $p_0$ has been removed recently by Zhang \cite{yz2018},
who proved 
that the limit $\nu$ is strictly positive if and only if $\PP(C=0)<1-\pc(\ZZ^d)$.
The problem of lower large-deviations of the $\Phi_{n,k}$ has been settled by Rossignol and Th\'eret \cite{rossth}
and Cerf and Th\'eret \cite{cerfth}.

For recent accounts of percolation and first-passage percolation, the reader is referred to 
\cite{MR3014795,auff-dam-han}.

\subsection{Word percolation}\label{sec:word}

The \lq word percolation' problem originated in work of Dekking \cite{dekk} on trees,
and was developed by Benjamini and Kesten as a generalisation of site percolation and
its variant called AB percolation. We start with the site percolation model on a given graph $G=(V,E)$
with density $p$, that is, a family $(X_v : v\in V)$ of independent Bernoulli variables
with $\PP_p(X_v=1)=1-\PP_p(X_v=0)=p$. Consider an alphabet of two letters $0$, $1$, 
and write $\Xi=\{0,1\}^\NN$
for the set of infinite \emph{words} in this alphabet.
A word $w=(w_i: i \in \NN) \in\Xi$ is said  to be \lq seen from vertex $v$'
if there exists a self-avoiding  walk $(v_1,v_2,v_3,\dots)$ from $v_1=v$ such that $X_{v_i}=w_i$ for all $i$;
$w$ is said to be \lq seen' if it is seen from some $v$.
The site percolation problem corresponds to the word $(1,1,1,\dots)$, and AB percolation to the word
$(1,0,1,0,\dots)$.

The questions of interest are: which words are seen from a given $v$ with strictly positive probability, and 
which words are seen somewhere in $G$ with probability $1$?
Harry wrote three papers on this topic with subsets of Benjamini and Sidoravicius/Zhang.
Let $S_v$ be the set of words seen from $v$, and let $S_\oo=\bigcup_v S_v$.

\begin{theorem}\label{thm:wordp}
Consider word percolation on the graph $G$. 
\begin{letlist}
\item {\bf\cite{MR1349161}}
Let $G=\ZZ^d$ and $\PP:=\PP_{\frac12}$. Then 
\begin{alignat*}{2}
\PP(S_\oo=\Xi)&=1 &&\q\text{if $d \ge 10$},\\
\PP(S_v=\Xi\text{ for some $v$})&=1 &&\q\text{if $d \ge 40$}.
\end{alignat*}

\item {\bf\cite{MR1637089}}
Let $G$ be the triangular lattice $\TT$ and $\PP:=\PP_{\frac12}$. With $\PP$-probability $1$,
every periodic word except $(1,1,1,\dots)$ and $(0,0,0,\dots)$, and almost
every non-periodic word, is seen. Here, \lq almost every word' means with $\mu_\beta$-probability $1$
where $\mu_\beta$ is product measure on $\Xi$ with some given density $\beta\in(0,1)$.

\item {\bf\cite{MR1825711}}
Let $G$ be obtained from $\ZZ^2$ by adding two diagonals to each face. 
We have $\PP_p(S_\oo=\Xi)=1$ whenever
$p\in(1-\pc,\pc)$, with $\pc$ $(<\frac12)$ the critical
probability of site percolation on $G$.

\end{letlist}
\end{theorem}

The word problem continues to fascinate probabilists, as exemplified by \cite{NTT},
where Theorem \ref{thm:wordp}(a) is extended to the statement that $\PP(S_\oo=\Xi)=1$ for $d \ge 3$.
Mention is made of three related
embedding problems which are proving to be challenging: the Lipschitz embeddings of \cite{grimhA,grimhB},
the quasi-isometric embeddings of \cite{basuss,basus}, and the Winkler compatibility problem \cite{MR3193962}. 

\subsection{Related processes}\label{sec:rel}

The percolation model is only the beginning of a rich modern theory of disordered discrete systems,
often indexed by lattices such as $\ZZ^d$. Brief accounts of Harry's contributions to 
a selection of related areas follow.

\subsubsection*{Greedy lattice animals}
To the vertices $v$ of $\ZZ^d$ are assigned iid
random variables $X_v$ satisfying $X_v \ge 0$. The \emph{weight} of a subset $U\subseteq V$ is
defined as the sum $W(U)=\sum_{v\in U} X_v$. Let $M_n$ be the maximal weight of connected subgraphs 
(also known as \lq lattice animals')  
with $n$ vertices and containing the origin. It was shown by Cox, Gandolfi, Griffin, and Kesten in 
\cite{MR1241039,MR1258174}
that, subject to a suitable moment condition on the $X_v$, the ratio $M_n/n$ converges as $n\to\oo$.
A similar result holds for $n$-step self-avoiding walks (SAWs) from the origin.  
Partial results without the positivity assumption on the $X_v$ are found in \cite{MR1825148}.  

\subsubsection*{$\rho$-percolation}
Consider site percolation on $\ZZ^d$ with $d \ge 2$ and density $p$, and let $\rho\in(0,1]$. 
Menshikov and Zuev \cite{menz}  initated the study of so-called $\rho$-percolation: 
we say that $\rho$-percolation occurs if there exists, with strictly positive probability, an infinite SAW 
$w=(w_0,w_1,\dots)$ from the origin $w_0=0$ such that the proportion $\pi_n$ of open vertices 
in the first $n$ vertices of $w$ satisfies $\liminf_{n\to\oo} \pi_n \ge \rho$. Thus, $\rho$-percolation
is connected to the greedy problem for SAWs above.

In \cite{MR1774071}, Kesten and Su studied $\rho$-percolation on the directed graph obtained from 
$\ZZ^d$ by orienting every edge according to increasing coordinate direction. Let $\pc(\rho)$
be the critical value of $p$ for given $\rho$.  They proved that
$$
d^{1/\rho}\pc(\rho) \to \frac{\theta^{1/\rho}}{e^\theta-1}\qq\text{as } d \to\oo,
$$
where $\theta=\theta(\rho)$ is the unique solution of a certain given equation.
This further example of Harry's interest in the large-$d$ behaviour of
interacting systems is complemented by  \cite{MR1002833} on the Ising model,
\cite{MR1070088} on Potts and Heisenberg models,
\cite{MR1064563} on percolation, and \cite{MR1144092} on the random-cluster model.

\begin{figure}
\centerline{\includegraphics[width=0.5\textwidth]{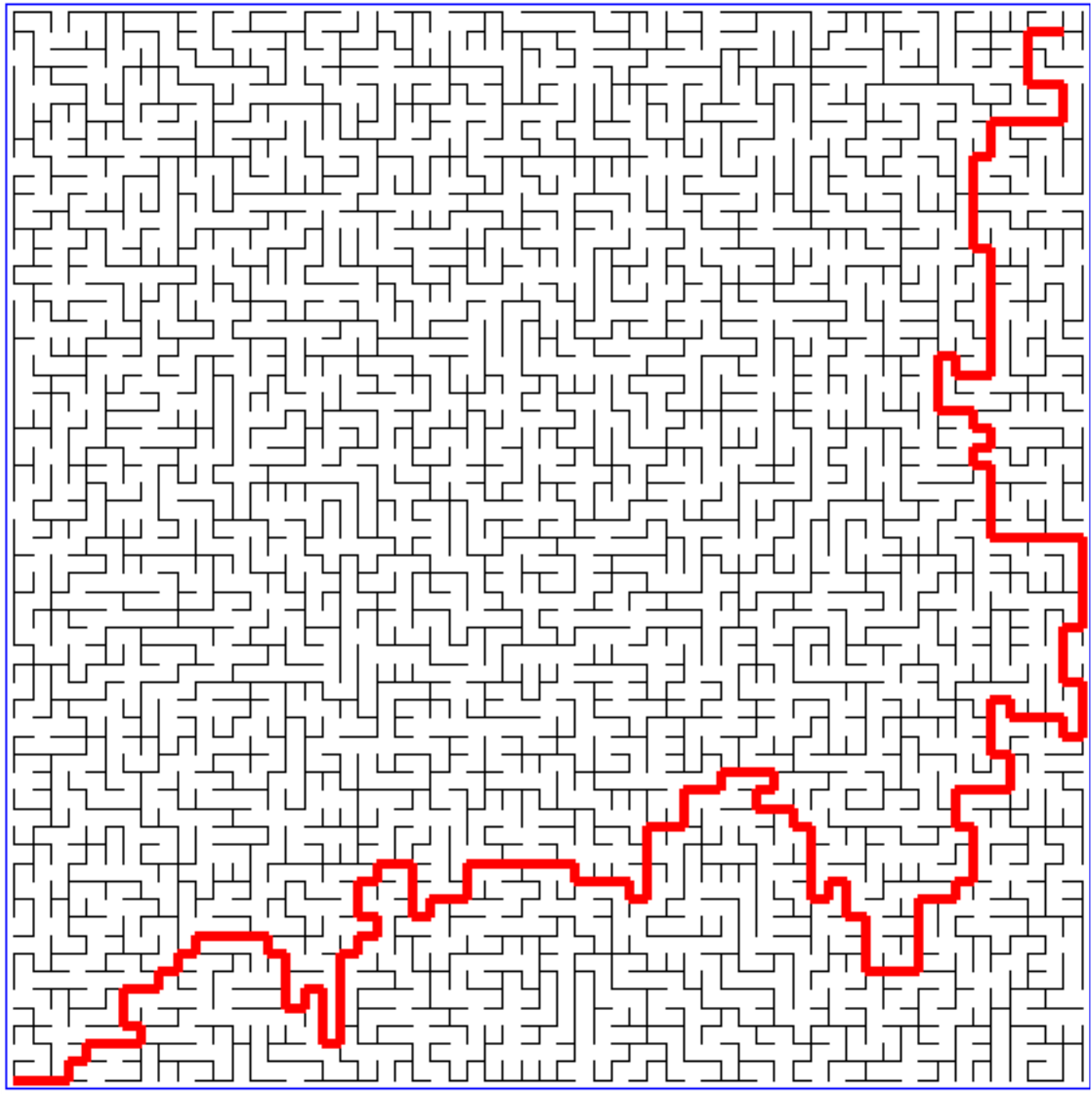}}
\caption{A uniform spanning tree (UST) on a square of
the square lattice. It contains a unique path between any two 
given boundary vertices, and this path has the law
of a loop-erased random walk. (By courtesy of Oded Schramm.)}
\label{fig:ust}
\end{figure}

\subsubsection*{Uniform spanning forest}
A \emph{uniform spanning tree} (UST) in a finite connected graph is a spanning tree chosen
uniformly at random from the set of all such trees. A \emph{uniform spanning forest} (USF) 
in $\ZZ^d$ is defined as the weak limit of USTs on larger and larger boxes
(see Figure \ref{fig:ust}). 
Pemantle \cite{pem} proved that USF is a.s.\ connected if and only if $d \le 4$.
The surprising fact was discovered by Benjamini, Kesten, Peres, and Schramm \cite{MR2123930}
that further transitions take place as $d$ passes through the multiples of $4$.
Suppose each tree in the USF is shrunk to a single vertex, and consider the graph $C$ thus obtained.
They showed that the diameter of $C$ increases by $1$ every time the dimension increases by $4$.
This unusual conclusion has been extended further by Hutchcroft and Peres \cite{hutchp}
in their proof that the USF with $d \ge 4$ undergoes a qualitative change to its connectivity 
\emph{every time} the dimension increases.

\subsubsection*{Minimal spanning tree}
Drop $n$ points $X_1,X_2,\dots,X_n$ uniformly at random into the unit cube $[0,1]^d$, where $d \ge 2$. 
Beardwood, Halton, and Hammersley \cite{BHH} initiated the study of graphical objects with vertex set $\{X_i\}$,
in their study (using subadditivity) of the length of the shortest travelling salesman path. 
Aldous and Steele \cite{aldste,ste88} asked about the minimal spanning tree on $\{X_i\}$, and proved that
its $\alpha$-length $M_{n,\alpha}$ (with an edge of length $l$ contributing length $l^\alpha$ for $\alpha\in(0,2]$)
is asymptotic to $n^{(d-\alpha)/d}$. 
Kesten and Lee \cite{MR1398055} proved a central limit theorem for this $M_{n,\alpha}$ 
by expressing it as a sum of martingale differences. Bounds on the rate of convergence in this last limit
were established recently by Chatterjee and Sen \cite{ChSen} using Stein's method.

\section{Further work}\label{sec:last}

\subsection{Probabilistic Diophantine approximation}\label{sec:eucl}

Harry wrote several works related to Diophantine approximation, frequently 
using probabilistic notation and methodology.
His first and only paper \cite{MR0097114} with Kac addressed a central limit theorem for the number of
appearances of a specified digit in the continued fraction expansion of a \lq typical' number.
This was achieved in the general context of the partial sums of a rapidly mixing sequence of iterates of a 
measure-preserving map.   The authors noted in an addendum that their results for continued fractions had
been obtained earlier by Doeblin \cite{doeb} in 1940.  

Let  $x,y \in [0,1]$, and let $0\le a\le b\le 1$. Define 
$$
S_n(x,y)= \sum_{k=1}^n 1\bigl(\{y+kx\}\in[a,b]\bigr),
$$
where $1(A)$ is the indicator function of $A$ and $\{z\}=z-\lfloor z\rfloor$ is the fractional part of $z$.
Weyl \cite{weyl} proved that $S_n(x,0)/n \to b-a$ as $n\to\oo$, for irrational $x$.
In \cite{MR113864,MR0142532}, 
Harry established a fluctuation theory of $S_n(X,Y)$, where $X$ and $Y$ are uniformly distributed on
$[0,1]$, namely that  
$$
\frac1{\log n} [S_n(X,Y)-n(b-a)] \Rightarrow Z/\rho,
$$ 
where $Z$ has the Cauchy distribution, and
 $\rho$  depends on $b-a$ and is constant for irrational $b-a$. This work was extended in \cite{MR168546}
 to an iterated logarithm law.
It is complemented in \cite{MR0210663,MR209253}
by a proof of a conjecture of Erd\H os and Sz\"usz that 
$S_n(x,0)-n(b-a)$ is bounded
in $n$ if and only if $b-a=\{jx\}$ for some integer $j$. This last result 
has motivated the
study of so-called bounded remainder sets (see \cite{DrTi}).

In a related work \cite{MR138612}, Ciesielski and Kesten  considered the process
$$
X_n(t) = \frac1{\sqrt n}\sum_{k=0}^{n-1} \bigl[ 1\bigl(\{2^kX\} \le t\bigr)-t\bigr], \qq t\in[0,1],
$$
where $X$ is uniformly distributed on $[0,1]$. 
They answered a question of Kac with their proof via an invariance principle that
$$
\lim_{n\to\oo} \PP\left(\sup_{t\in[0,1]} |X_n(t)| <u\right) 
=\PP\left(\sup_{t\in[0,1]} |X(t)|<u \right), \qq u\in [0,1],
$$
where $X$ is the Gaussian process that is the weak limit $\lim_{n\to\oo}X_n$.

Let $\langle x \rangle$ denote the distance from $x$ to its closest integer. In an extension \cite{MR137692} of work of
Friedman and Niven \cite{friniv} and Erd\H os, Sz\" usz, and Tur\' an \cite{erdos64,est}, 
Harry showed three theorems, including
that, if $X$ is uniformly distributed on $[0,1]$,
then $m \cdot\min\{\langle kX\rangle: 1\le k \le m\}$ converges in distribution as $m\to\oo$. 
Recent developments of Athreya and Ghosh may be found in \cite{AG18}.

\subsection{Random growth}

First-passage percolation is an example of a random growth model, and there are many others, including the
classic Eden model \cite{eden}. Kesten and Schonmann
\cite{MR1327220} explored a variant of the Eden model, proving an asymptotic growth rate and shape.

Kesten and Sidoravicius \cite{MR1961167,MR2184100,MR2247840,MR2440925,MR2415386} studied a population model
on $\ZZ^d$ with two types of particle. Type-A particles move as independent random walks in continuous time
with rate $D_A$, and
type-B particles behave similarly at rate $D_B$. When a type-B particle meets a type-A particle, the latter changes 
its type to B instantaneously. The problem was to prove a shape theorem for the spread of type-B particles.
This process has been called the \emph{frog model} in the special case $D_A=0$. 
Recent work on a continuum variant of 
this problem includes \cite{BDDHJ,GZL}.

\subsection{Population genetics}

Inspired by a lecture course given by Kingman in 1979, 
Harry dabbled briefly in population genetics, with work on the so-called Ohta--Kimura model
for the evolution of allelic frequencies in a certain finite (but large) population.
Each generation has a fixed size $N$, from which a random sample of size $n$ is taken. 
Let $\Lambda(n,N,t)$ be the number of distinguishable alleles present in such a sample from 
generation $t$, and let $\Lambda(n,N)$ be
the weak limit of $\Lambda(n,N,t)$ as $t \to\oo$  (this limit exists by a result of Kingman \cite{kingm}).
Kimura and Ohta's analysis suggested that $\EE\Lambda(N,N)$ remains bounded as $N\to\oo$,
whereas Moran showed heuristically that $\Lambda(N,N)\to\oo$.
It turned out that, in a fashion, both were correct.

Let $\gamma_k$ denote tetration (or the power tower) of height $k$, $\gamma_k := e^{e^{\cdot^{\cdot^{e}}}}$, and let
$\lambda(n)=\max\{k: \gamma_k\le n\}$ denote its inverse. Harry showed in \cite{MR602144,MR596465}
that, 
in the limit as $n,N\to\oo$ with $n \le N$, $\Lambda(n,N)/\lambda(n)$ converges in probability to
an explicit limit. He pointed out that, for all practical purposes, $\Lambda(N,N)$ remains bounded, 
and he illustrated this with the observation that, for $3814280 \le N \le 10^{1656520}$, we have $\lambda(N)=3$. 
After lecturing on this at Stanford University, he commented 
on the evident mystification of the biologists present, and on
Karlin's expression of satisfaction at having 
predicted the boundedness of $\Lambda(N,N)$.

\subsection{Quasi-stationary distributions of Markov chains}

Let $X$ be  a discrete-time Markov chain with an absorbing state labelled $0$. What can be said about $X_n$
conditional on $X_n\ne 0$? This classical question has led to a theory of so-called quasi-stationary distributions
(qsd).  Harry, in combinations with Ferrari, Mart{\'\i}nez, and Picco, proved several fundamental results
concerning qsds. They proved in \cite{MR1334159} that, subject to a certain condition,
 a qsd exists if and only if the time to absorption 
has an exponentially-decaying tail. The required condition fails for an important class of Markov chains
arising from spatially distributed interacting particle systems
including the (discrete-time) subcritical contact model on $\ZZ^d$. This was rectified in  the significant work
\cite{MR1398060}, with a study based around the property of so-called $R$-positive-recurrence.
On his way to these results, Harry \cite{MR1439524}
revisited the topic of ratio limit theorems, with proofs that ratio limits exist for
Markov chains on the non-negative integers with an absorbing state $0$, conditional
on not yet being absorbed.

\subsection{Annihilating and coalescing random walks}
Interacting random walks of these two types have been studied since the 1970s as duals to the
antivoter and voter models of interacting particle systems. In a model first studied by 
van den Berg and Kesten \cite{MR1756007}, particles are placed at each vertex of $\ZZ^d$, 
and each moves in the manner of a Markov chain. 
When a particle jumps to a vertex already occupied by $j$ particles, it
is removed with some given probability $q_j$. Let $p(t)$ be the mean density of occupied vertices at time $t$.
The authors showed that $p(t)\sim C/t$ if $d \ge 6$. Subject to an extra condition, 
this was extended in \cite{MR1901947} to $d \ge 3$. 

This project was continued by Harry in \cite{MR1771527} with a proof that $p(t)\sim C/t$
when $d\ge 9$ for similar models in which the annihilation/coalescence can take place within the
set of neighbours of the moving particle.

In one of Harry's final papers, he returned with  Benjamini, Foxall, Gurel-Gurevich, and Junge \cite{MR3522593}
to coalescing random walks, this time in the context of
a connected, locally finite graph $G$. They asked when  such walks are \lq site recurrent'
in that every site is a.s.\ visited infinitely often. Site recurrence is equivalent to $\int_0^\oo p_v(t)\,dt=\oo$
for all vertices $v$, where $p_v(t)$ is the probability that $v$ is occupied at time $t$.
They showed that $p_v(t) \ge C/(1+t)$ when $G$ has bounded degree, and $p_v(t)\ge C/(t \log t)$ when $G$ is a 
branching process whose family-sizes have an exponentially decaying tail, thereby verifying site recurrence in both cases. 

\section*{Acknowledgements}
\addcontentsline{toc}{section}{Acknowledgements}
The author thanks Jean Bertoin, Rick Durrett, and Ross Maller 
for their helpful comments on this article, and he gratefully acknowledges
the thoughtful and detailed report of an anonymous referee.



\newpage
\renewcommand\refname{Publications of Harry Kesten}

\end{document}